\documentclass[final,hidelinks,onefignum,onetabnum]{siamart250211}

\usepackage{enumerate}
\usepackage{bm}
\usepackage{amsmath}
\usepackage{amssymb}
\usepackage{mathrsfs}
\usepackage{graphicx}
\usepackage[version=4]{mhchem}
\usepackage{siunitx}
\usepackage{cleveref}
\usepackage{todonotes}

\ExplSyntaxOn
\keys_define:nn { mhchem }
{
	arrow-min-length .code:n =
	\cs_set:Npn \__mhchem_arrow_options_minLength:n { {#1} } 
}
\ExplSyntaxOff
\mhchemoptions{arrow-min-length=1em}

\newcommand{\R}{\mathbb{R}}
\newcommand{\Id}{\mathbb{I}}
\let\div\relax \DeclareMathOperator{\div}{div}
\DeclareMathOperator{\spn}{span}
\DeclareMathOperator{\tr}{tr}
\newcommand{\norm}[1]{\left\lVert#1\right\rVert}

\newsiamremark{assumption}{Assumption}
\crefname{assumption}{Assumption}{Assumptions}

\headers{FEM for electroneutral multicomponent flows}
{Aaron Baier-Reinio, Patrick E.~Farrell and Charles W.~Monroe}

\title{Finite element methods for electroneutral multicomponent electrolyte flows%
\thanks{\funding{ABR was supported by a Clarendon scholarship from 
the University of Oxford. PEF was supported by EPSRC grants 
EP/R029423/1 and EP/W026163/1, and by the Donatio Universitatis
Carolinae Chair ``Mathematical modelling of multicomponent systems''.
CWM was supported by the Faraday Institution Multiscale Modelling Project, grant FIRG059.}}}

\author{Aaron Baier-Reinio\thanks{Mathematical Institute, University of Oxford, 
	Oxford, OX2 6GG, UK
	(\email{aaron.baier-reinio@maths. ox.ac.uk}).}
	\and Patrick E.~Farrell\thanks{Mathematical Institute, University of Oxford, 
	Oxford, OX2 6GG, UK and
    Mathematical Institute, Faculty of Mathematics and Physics, Charles 
    University, Czechia
    (\email{patrick.farrell@maths.ox.ac.uk}).}
	\and Charles W.~Monroe\thanks{Department of Engineering Science, University of 
	Oxford, Oxford, OX1 3PJ, UK and The Faraday
	Institution, Harwell Campus, Didcot, OX11 ORA, UK
	(\email{charles.monroe@eng.ox.ac.uk}).}}

\ifpdf
\hypersetup{
  pdftitle={Finite element methods for electroneutral multicomponent flows},
  pdfauthor={Aaron Baier-Reinio, Patrick E.~Farrell, and Charles W.~Monroe}
}
\fi

\begin{document}

\maketitle

\begin{abstract}
We present a broad family of high-order finite element algorithms for simulating 
the flow of electroneutral electrolytes.
The governing partial differential equations that we solve are the
electroneutral Navier--Stokes--Onsager--Stefan--Maxwell (NSOSM) equations,
which model momentum transport, multicomponent diffusion and electrical effects
within the electrolyte.
Our algorithms can be applied in the steady and transient settings, in two and 
three spatial dimensions, and under a variety of boundary conditions.
Moreover, we allow for the material parameters (e.g.~viscosity, diffusivities, 
thermodynamic factors and density) to be dependent on the local state of the 
mixture and thermodynamically non-ideal.
The flexibility of our approach requires us to address subtleties that arise 
in the governing equations due to the interplay between boundary conditions and 
the equation of state.
We demonstrate the algorithms in various physical configurations,
including (i) electrolyte flow around a microfluidic rotating disk electrode
and (ii) the flow in a Hull cell of a cosolvent electrolyte mixture used
in lithium-ion batteries.
\end{abstract}



\section{Introduction} \label{sec:introduction}

We address the numerical simulation of liquid electrolytes, which are fluids
that transport electrical charge by the motion of ions.
Electrolytes are \textit{multicomponent fluids} (or \textit{mixtures}),
which means that they consist of multiple distinct chemical species
\cite{giovangigli2012multicomponent,wesselingh2000mass}.
While a general mixture may consist entirely of uncharged species,
in an electrolyte at least two of the species must be ions of opposite charge.
For example, dissolving table salt in water yields an electrolyte with
three species: \ce{H2O}, \ce{Na+} and \ce{Cl-}.
Prominent applications of electrolytic flows include energy storage
(e.g.~batteries and fuel cells),
chemical processes (e.g.~electrodialysis and electroplating)
and biological systems (e.g.~biological membranes)
\cite{bartlett2008bioelectrochemistry,newman2021electrochemical,wesselingh2000mass}.

We assume that the electrolyte is \textit{electroneutral}, which means that
its local charge density is everywhere zero%
\footnote{Note that electroneutrality does not prevent the transport of charge
within the mixture.}.
This is a common and accurate assumption for most
electrochemical systems at length scales much larger than the Debye length
(which is typically on the order of nanometers)
\cite{newman2021electrochemical,wesselingh2000mass}.
For example, in lithium-ion batteries, electroneutrality holds throughout the bulk
of the electrolyte and is only violated in sub-nanometer-wide layers near the
electrolyte-electrode interfaces.
Consequentially, physics-based microscale lithium-ion battery models typically use
differential equations that assume electroneutrality throughout the
domain and incorporate interfacial charge layers through boundary conditions
\cite{brosa2022continuum,richardson2022charge}.
In this work, electroneutrality allows us to follow work done earlier by Van-Brunt et al.\ \cite{van2023structural}
and employ a change of basis that transforms the governing equations into a form that is structurally similar to
that of an uncharged mixture, simplifying both the model and numerics.

A constitutive relation for the
species mass fluxes must be chosen to model mass transport in multicomponent flows.
For electrolytic flows, the popular Nernst--Planck model accounts for mass transport by convection, Fickian diffusion and 
electromigration \cite{bard2022electrochemical,dreyer2013overcoming}.
Several finite element methods exist for the electroneutral Nernst--Planck
model \cite{bortels1996multi,ellingsrud2025splitting,roy2023scalable},
and literature for its non-electroneutral analogue, the
Poisson--Nernst--Planck model, is especially abundant
(see e.g.~\cite{correa2023new,prohl2010convergent,xie2020effective}).
However, because the Nernst--Planck model assumes that every species in isolation obeys Fick's law of diffusion
\cite{fick1855ueber}, it has the drawback of being unable to capture
\textit{cross-diffusion}, a physical phenomenon that arises
when different species exert diffusional forces on each other
\cite{jungel2015boundedness,sun2019entropy,vanag2009cross}.
For \textit{dilute mixtures}, i.e., mixtures
with only one species present in non-trace amounts,
it is usually assumed appropriate to neglect cross-diffusion, and doing so
greatly simplifies the model.
However, many practical problems involve non-dilute mixtures for which
Fick's law is inadequate.
For example, in lithium-ion battery electrolyte modelling, cation-anion 
cross-diffusion must be accounted for to match experimentally observed 
conductivity, salt diffusivity and transference numbers 
\cite{bizeray2016resolving}.
For detailed discussions on the limitations of Fick's law, we refer to 
\cite{krishna2015uphill,krishna2019diffusing,krishna1997maxwell}
for generic mixtures and
\cite[Chapters 11-12]{newman2021electrochemical} for electrolytes.
Drawbacks of the Nernst--Plack model from a 
thermodynamical perspective are also given in \cite{dreyer2013overcoming}.

To avoid the limitations of Nernst--Planck theory,
we treat mass transport using the Onsager--Stefan--Maxwell%
\footnote{Nomenclature varies depending on the source; sometimes the name
\textit{Maxwell--Stefan} or even
\textit{generalized Maxwell--Stefan} is used instead.}
(OSM) constitutive laws \cite{bird2002transport,krishna1997maxwell,wesselingh2000mass}.
These equations are based on irreversible thermodynamics
\cite{de2013non,onsager1931reciprocal,onsager1931reciprocal2}
and generalize the Stefan--Maxwell equations, which model cross-diffusion in
ideal gases \cite{maxwell1866,stefan1871gleichgewicht}.
In electrochemistry, the OSM equations were popularized by Newman
\cite{newman2021electrochemical,newman1965mass}
for modelling electrolytes, and the resulting 
framework is sometimes known as \textit{concentrated solution theory}.
The OSM equations provide a thermodynamically rigorous model for mass 
transport by convection, cross-diffusion and electromigration.
The equations are also compatible with electroneutrality \cite{van2023structural}
and anisothermality \cite{van2022consolidated},
although we do not treat the latter here.
Moreover, the OSM equations can handle thermodynamic
non-idealities (arising in mixtures with a nonzero excess Gibbs free energy
of mixing \cite{doble2007perry}), the effects of which are important
in the non-dilute regime \cite{krishna1997maxwell}.
In settings where momentum transport is important
(e.g.~fuel cells, electrodialysis),
it is necessary to couple the OSM equations to momentum
conservation laws such as the Stokes or Navier--Stokes equations; we respectively call
these coupled sets of constitutive laws the SOSM or NSOSM equations.

Due to its flexibility in treating complicated transport phenomena and nonideal
thermodynamics, the OSM framework has received much attention in
electrochemistry as a tool for modelling electrolytic flows
\cite{newman2021electrochemical}.
For example, the framework has been applied to lithium-ion batteries
\cite{richardson2022charge},
fuel cells \cite{weber2014critical,weber2008mathematical},
electrodialysis \cite{kraaijeveld1995modelling}
and electrolysis \cite{sijabat2019maxwell}.
Nonetheless, the OSM equations have received limited attention in the scientific
computing literature, especially for the coupling of OSM diffusion to Stokes or
Navier--Stokes flow.
Most numerics papers on the OSM equations assume a thermodynamically ideal
gaseous mixture, such as the finite element methods in
\cite{braukhoff2022entropy,carnes2008local,jungel2019convergence,mcleod2014mixed,van2022augmented}.
For numerics papers on the SOSM or NSOSM equations, we are aware of
\cite{ern1998thermal} for a finite difference scheme,
\cite{balakrishnan2014fluctuating,bhattacharjee2015fluctuating,donev2014low,donev2015low,
donev2019fluctuating,peraud2016low,srivastava2023staggered}
for a series of Cartesian-grid finite volume schemes under various physical 
assumptions in the setting of fluctuating hydrodynamics, and
\cite{aznaran2024finite,baier2024high,brunk2025structure,burman2003bunsen,longo2012finite}
for finite element schemes.
Among these, only \cite{baier2024high} considers spatially high-order methods, and 
only \cite{donev2019fluctuating,peraud2016low}
consider charged species 
(although \cite{donev2019fluctuating,peraud2016low} make the Boussinesq 
approximation to simplify the numerics).
We are not aware of existing finite element algorithms for the electroneutral 
NSOSM equations, a gap the present work addresses.

A major novelty of this work lies in our treatment of the electroneutrality 
condition, which in mathematical terms places an algebraic constraint on the species 
concentrations.
Previous works have used this constraint to obtain an elliptic equation for 
the electrical potential, which is then coupled to the relevant equations that 
govern mass and momentum transport
\cite{donev2019fluctuating,ellingsrud2025splitting,roy2023scalable}.
The structure of these coupled transport equations does not resemble that afforded by an uncharged mixture, requiring the development of new solution techniques:
temporal splitting schemes are introduced in
\cite{donev2019fluctuating,ellingsrud2025splitting}
and a monolithic approach is considered in \cite{roy2023scalable}.
By contrast, in this work we use the so-called \textit{salt-charge transformation},
introduced in \cite{van2023structural} and inspired by Newman's work on binary 
electrolytes \cite{newman1965mass}, to transform the electroneutral NSOSM
equations into a form that structurally resembles the NSOSM equations for an
uncharged mixture.
We then apply the method of lines, whereby we spatially discretize the
transformed problem using a steady NSOSM solver and then solve the resulting
spatially discrete system of differential-algebraic equations by time-stepping.
Hence, our approach has the advantage of allowing uncharged NSOSM spatial
discretizations and solvers to be employed for electroneutral electrolytic mixtures.

The steady NSOSM discretization we employ modifies the SOSM scheme introduced in \cite{baier2024high}, now including the 
nonlinear convective terms in the momentum balance.
These are trivial to include because this work assumes low Reynolds number flow, so 
the convective terms do not require stabilization.
We also apply a second modification that is straightforward but has 
consequences regarding the set of material parameters one must specify to implement the model.
Namely, instead of discretizing the OSM equations using
chemical potentials (as is done in 
\cite{aznaran2024finite,baier2024high,brunk2025structure}),
we express the OSM diffusional driving forces in terms of mole fraction gradients
and we only discretize the mole fractions.
Because we do not discretize the chemical potentials, we do not require 
experimental knowledge of species activity coefficients in the mixture. Only the 
\textit{thermodynamic factors} (i.e.~partial derivatives 
of activity coefficients with respect to mole fractions) are required,
which allows us to straightforwardly draw on experimental data from
\cite{wang2021potentiometric,wang2020shifting}
in our numerical experiments.

When modelling mixtures, the interplay between the volumetric 
equation of state (EOS) and boundary conditions (BCs) is subtle.
The EOS and BCs can both impose restrictions on the species
concentrations, and these restrictions must be compatible.
Investigating this interplay is mathematically challenging even for
a two-species Fickian mixture \cite{feireisl2016mathematical}.
Here, this challenge is exacerbated by the fact that our steady
NSOSM solver requires user-chosen integral constraints to be imposed
on the mathematical solution%
\footnote{
Henceforth we shall always use ``solution'' in the mathematical, rather than the physical, sense.
}%
\cite{baier2024high}.
These constraints arise naturally in the steady case, because the steady
mass continuity equations do not dictate how many moles of each species are
present in the domain; this must be prescribed by additional integral constraints.
However, in the transient case, the constraints must be chosen carefully, such that
they are compatible with the choice of EOS and BCs.
We discuss several physically relevant choices of EOS and BCs
along with appropriate corresponding constraints. We also discuss a
seemingly physical combination of EOS and BCs that appears to yield an ill-posed
problem. We are not aware of discussions regarding this interplay in the multicomponent
numerics literature and believe these findings may be of independent interest.
Incidentally, we mention that a discussion about the challenges of
discretizing the multicomponent mass continuity equations, in a way that is
compatible with the EOS, is given in \cite{donev2014low} in the case of a
low-order finite volume scheme%
\footnote{The challenges encountered in \cite{donev2014low} are different from 
those encountered in this work, presumably because they discretize the species
partial densities whereas we discretize mole fractions, and also because they 
attempt to enforce the EOS pointwise whereas we do not.}.

The remainder of this paper is organized as follows.
In \cref{sec:cts} we introduce the electroneutral NSOSM equations, and show how the
salt-charge transformation can be applied to obtain an equivalent system
amenable to a steady NSOSM solver.
In \cref{sec:spatial_disc} we first consider the transformed steady problem, which
we discretize using the (modified) scheme of \cite{baier2024high}.
In \cref{sec:temporal_disc} we apply the method of lines to the transformed
transient problem, discretizing in space using the steady scheme from
\cref{sec:spatial_disc} and taking special care to identify what
constraints must be imposed depending on the EOS and BCs.
In \cref{sec:numerical} we demonstrate our algorithm with a variety of numerical 
examples, and conclusions are drawn in \cref{sec:conc}.

\section{Governing equations and salt-charge transformation} \label{sec:cts}

We consider an isothermal, chemically nonreacting mixture of 
$n \geq 3$ species indexed by $i \in \{1 : n\}$.
Let $z_i$ denote the equivalent charge of species $i$ and
$n_c \geq 2$ the number of charged species.
Note that $n_c=n$ is possible, e.g.~for a molten salt.
We list the species so that the first $(n-n_c)$ species
are uncharged and the last two species are oppositely charged.
This means that $z_1, \ldots, z_{n-n_c} = 0$ and
$z_{n-n_c+1}, \ldots, z_n \neq 0$ with $z_n / z_{n-1} < 0$.

We next state the governing partial differential equations, which are posed
on a bounded and connected Lipschitz spatial domain
$\Omega \subset \R^d$ with $d \in \{2, 3\}$.

\subsection{Electroneutral NSOSM equations} \label{sec:cts_electroneutral}

The thermodynamic state variables that we consider are the temperature $T$,
the pressure $p$ and the species molar concentrations $c_i$ for
$i \in \{1 : n \}$.
Since the mixture is isothermal $T$ is a fixed parameter,
but $p$ and $c_1, \ldots, c_n$ may vary with space and time.
The total concentration is $c_T = \sum_{j=1}^n c_j$
and the density is $\rho = \sum_{j=1}^n m_j c_j$ where $m_i$ is the molar mass of
species $i$.
Composition can be described using the species mole fractions
$x_i = c_i / c_T$.
Note that by definition $\sum_{j=1}^n x_j = 1$.
The charge density is
$\rho_e = F \sum_{j=1}^n z_j c_j$ where $F$ is Faraday's constant.
Electroneutrality requires that $\rho_e = 0$, or equivalently
\begin{equation} \label{eq:ensosm_en}
	\textstyle\sum_{j=1}^n z_j c_j = 0 \quad \textrm{in } \Omega,
\end{equation}
which restricts the space of permissible compositions.

To model momentum and mass transport we employ the Cauchy momentum
equation, and nonreacting species mass continuity equations,
\begin{alignat}{2}
	\frac{\partial (\rho v)}{\partial t}
	+ \nabla \cdot (\rho v \otimes v)
	+ \nabla p
	- \nabla \cdot \tau &= \rho f
	\quad &&\textrm{in } \Omega, \label{eq:ensosm_ns} \\
	\frac{\partial c_i}{\partial t}
	+ \nabla \cdot N_i &= 0
	\quad &&\textrm{in } \Omega
	\quad \forall i \in \{1 : n\}. \label{eq:ensosm_mc}
\end{alignat}
Here $v$ is the barycentric velocity,
$\tau$ is the viscous stress, $f$ is a prescribed body force
which is typically zero,
and $N_i$ is the molar flux of species $i \in \{1 : n\}$.
The barycentric velocity is related to the molar fluxes through the
\textit{mass-average constraint}
\begin{equation} \label{eq:ensosm_mavg}
	\rho v = \textstyle\sum_{j=1}^n m_j N_j \quad \textrm{in } \Omega.
\end{equation}
Since $\rho = \sum_{j=1}^n m_j c_j$ one can multiply \cref{eq:ensosm_mc}
by $m_i$ and sum over $i$ while using \cref{eq:ensosm_mavg} to obtain the equation
for total mass continuity
$\partial_t \rho + \nabla \cdot (\rho v) = 0$.
However, we never explicitly discretize this equation since it is a consequence
of \cref{eq:ensosm_mc} and \cref{eq:ensosm_mavg}.

The transport equations 
\cref{eq:ensosm_ns,eq:ensosm_mc}
must be closed with constitutive laws for the
viscous stress and mass fluxes.
For the former we apply the Newtonian fluid law
\begin{equation} \label{eq:ensosm_newtonian}
	\tau = 2 \eta \epsilon(v)
	+ (\zeta - 2 \eta / d) (\nabla \cdot v) \Id,
\end{equation}
where $\epsilon(v)$ is the symmetric gradient of $v$,
$\Id$ is the $d \times d$ identity matrix,
$\eta > 0$ is the shear viscosity and $\zeta > 0$ is the bulk viscosity.
We allow the viscosities $\eta$ and $\zeta$ to depend on the thermodynamic
state variables $T, p$ and $x_1, \ldots, x_n$.

To model the mass fluxes we employ the isothermal, non-isobaric
Onsager--Stefan--Maxwell (OSM) equations 
\cite{bird2002transport,krishna1997maxwell,newman2021electrochemical,wesselingh2000mass},
which can be written using molar fluxes as
\begin{equation} \label{eq:ensosm_osm}
	-\nabla \mu_i + \frac{m_i}{\rho} \nabla p = 
	\textstyle\sum_{j=1}^n \bm{M}_{ij} N_j
	\quad \textrm{in } \Omega
	\quad \forall i \in \{1 : n\},
\end{equation}
where $\mu_i$ is the \textit{electrochemical potential} of species $i \in \{1:n\}$
and $\bm{M}$ is the \textit{Onsager transport matrix}.
Two important properties of $\bm{M}$ are that it is symmetric positive-semidefinite,
and its nullspace is spanned by $[c_1, \ldots, c_n]^{\top}$.
Moreover, the Gibbs--Duhem relation implies that
$\sum_{j=1}^n c_j (-\nabla \mu_j + [m_j / \rho]\nabla p) = 0$
\cite{guggenheim1967thermodynamics}.
The left side of \cref{eq:ensosm_osm} therefore lies in the
orthogonal complement of $\spn\{[c_1, \ldots, c_n]^{\top}\}$,
which is the range of $\bm{M}$.
Hence \cref{eq:ensosm_osm} is non-uniquely solvable for $N_1, \ldots, N_n$.
Uniqueness of the fluxes comes from imposing the mass-average constraint 
\cref{eq:ensosm_mavg} in addition to \cref{eq:ensosm_osm}.
Often $\bm{M}$ is written in terms of so-called 
\textit{Stefan--Maxwell diffusivities}
$\mathscr{D}_{ij}$ as
\begin{equation} \label{eq:transport_m_def}
	\bm{M}_{ij} =
	\begin{cases}
		-\frac{RT}{\mathscr{D}_{ij} c_T}
		& \text{if $i \neq j$}, \\
		\textstyle\sum_{k=1, k \neq i}^n \frac{RT c_k}{\mathscr{D}_{ik} c_T c_j}
		& \text{if $i = j$},
	\end{cases}
\end{equation}
with $R$ the gas constant.
The diffusivities satisfy
$\mathscr{D}_{ij}=\mathscr{D}_{ji}$ for $i \neq j$ \cite{bird2002transport},
while $\mathscr{D}_{ii}$ is undefined. We allow the $\mathscr{D}_{ij}$ 
to depend on $T, p$ and $x_1, \ldots, x_n$
(see e.g.~\cite{kraaijeveld1993negative}).

\Crefrange{eq:ensosm_en}{eq:ensosm_osm} comprise the electroneutral
NSOSM equations.
In addition to needing suitable boundary and initial conditions, the equations
as written are still not closed.
One must additionally supply a volumetric
equation of state (EOS), which gives the total concentration $c_T$
(or equivalently, mass density $\rho$) as a function of $T, p$ and $x_1, \ldots, x_n$.
A constitutive law giving the electrochemical potentials $\mu_i$ as a function
of $T, p, x_1, \ldots, x_n$ and the 
\textit{electrical state} of the system must also be given.
However, in multicomponent systems, formulating the notion of
electrical state is delicate, and requires a foray into
technicalities of electrochemical thermodynamics.

Intuitively, the electrical state should quantify how much electrical potential
energy is locally available in the system.
A reasonable model may be to introduce an ``electrostatic potential''
$\phi$ which acts as a Lagrange multiplier for enforcing electroneutrality,
and to then decompose the electrochemical potentials as
$\mu_i = \mu_i^{\textrm{chem}} + \phi$, where $\mu_i^{\textrm{chem}}$
encodes the ``chemical'' or ``non-electrical'' part of $\mu_i$
and depends on $T, p$ and $x_1, \ldots, x_n$ only.
However, Guggenheim and Newman have criticized such decompositions as being
ambiguous \cite{guggenheim1967thermodynamics,newman2021electrochemical},
with Newman noting that several inequivalent notions of ``electrostatic
potential'' exist in concentrated mixtures
(see also \cite{boettcher2020potentially,pethica2007electrostatic}).
Fortunately, we do not need to contend with this thermodynamical 
technicality in the present work.
Indeed, owing to electroneutrality, the salt-charge transformation will eliminate
any need to model the dependence of electrochemical potentials on the electrical
state \cite{van2023structural}.
In what follows, it will only be necessary to model the dependence of 
\textit{chemical potentials} of hypothetical \textit{uncharged salts}
on $T, p$ and $x_1, \ldots, x_n$%
\footnote{In the non-electroneutral setting this strategy does not work. 
One must instead choose how to quantify the electrical state and model how the
electrochemical potentials depend on it.}.

\subsection{Salt-charge transformation} \label{sec:cts_saltcharge}

It will be helpful to use boldface notation for vectors and matrices that
are indexed by $i \in \{1 : n\}$.
For example, we shall write
$\bm{z} = [z_1, \ldots, z_n]^{\top}$,
$\bm{c} = [c_1, \ldots, c_n]^{\top}$,
$\bm{m} = [m_1, \ldots, m_n]^{\top}$
and so on.
These should be interpreted as column vectors (or $n \times 1$ matrices).
Taking their gradient yields $n \times d$ matrices, hence
for example $(\nabla \bm{\mu})_{ij} = (\nabla \mu_i)_j$.
For the fluxes and other $\R^d$-valued quantities we write
$\bm{N} = [N_1, \ldots, N_n]^{\top}$ which is understood to be an
$n \times d$ matrix, and the divergence acts on its rows so that
$\nabla \cdot \bm{N}$ is a column vector (or $n \times 1$ matrix).
With this notation, the governing equations in
\cref{eq:ensosm_en,eq:ensosm_mc,eq:ensosm_mavg,eq:ensosm_osm} become
\begin{alignat}{2}
	\bm{z}^{\top} \bm{c} &= 0
	\quad &&\textrm{in } \Omega,
	\label{eq:v_ensosm_en} \\
	\frac{\partial \bm{c}}{\partial t}
	+ \nabla \cdot \bm{N} &= 0
	\quad &&\textrm{in } \Omega,
	\label{eq:v_ensosm_mc} \\
	v &= \bm{\psi}^{\top} \bm{N}
	\quad &&\textrm{in } \Omega,
	\label{eq:v_ensosm_mavg} \\
	-\nabla \bm{\mu} + \bm{\psi} \nabla p &= \bm{M} \bm{N}
	\quad &&\textrm{in } \Omega,
	\label{eq:v_ensosm_osm}
\end{alignat}
where $\bm{\psi} := \bm{m} / \rho$ and
$\nabla p$ is viewed as being a $1 \times d$ row vector.

We now carry out the salt-charge transformation.
While a detailed presentation is given in \cite{van2023structural},
we briefly outline it here for completeness,
and because \cite{van2023structural} assumes an isobaric mixture and hence does 
not keep
track of pressure gradient terms in the OSM equations.
The starting point is an $n \times n$ transformation matrix $\bm{Z}$ given by
\begin{equation} \label{eq:sc_trans_mat}
	\bm{Z} = \begin{bmatrix}
		\bm{\nu}_1^{\top} \\
		\vdots \\
		\bm{\nu}_{n-1}^{\top} \\
		\bm{z}^{\top} / \norm{\bm{z}}
	\end{bmatrix}.
\end{equation}
Here $\norm{\bm{z}} = \sqrt{\bm{z}^{\top} \bm{z}}$
and $\bm{\nu}_1, \ldots \bm{\nu}_{n-1}$ are $n \times 1$ column vectors
whose entries contain the stoichiometric coefficients of 
$n-1$ independent hypothetical chemical reactions that neutralize the species.
An example from \cite{van2023structural} is an $n=5$ species mixture
consisting of \ce{H2O}, \ce{Na+}, \ce{Cl-}, \ce{Mg^{2+}} and \ce{SO^{2-}_4}.
Therefore $\bm{z}^{\top} = [0, 1, -1, 2, -2]$, and
a choice of neutralizing reactions with corresponding stoichiometric coefficient
vectors is
\begin{alignat*}{2}
	\ce{H2O &<=> H2O}
	\qquad && \longrightarrow \qquad \bm{\nu}_1^{\top} = [1, 0, 0, 0, 0], \\
	\ce{Na+ + Cl- &<=> NaCl}
	\qquad && \longrightarrow \qquad \bm{\nu}_2^{\top} = [0, 1, 1, 0, 0], \\
	\ce{Mg^{2+} + 2Cl- &<=> MgCl_2}
	\qquad && \longrightarrow \qquad \bm{\nu}_3^{\top} = [0, 0, 2, 1, 0], \\
	\ce{2Na+ + SO^{2-}_4 &<=> Na_2SO_4}
	\qquad && \longrightarrow \qquad \bm{\nu}_4^{\top} = [0, 2, 0, 0, 1].
\end{alignat*}
Note that the choice of reactions may not be unique.
Following \cite{van2023structural} we make the convention that
$\bm{\nu}_1, \ldots \bm{\nu}_{n-1}$ have the following properties.
First, for the uncharged species $i \in \{1:n-n_c\}$ we have
$(\bm{\nu}_i)_j = \delta_{ij}$, corresponding to a trivial identity reaction.
Second, for $i \in \{n-n_c+1:n-1\}$ the first $n-n_c$ entries of $\bm{\nu}_i$
are zero, and exactly two of the remaining entries
of $\bm{\nu}_i$ are nonzero, and must be coprime positive integers such that
$\bm{\nu}_i$ is orthogonal to $\bm{z}$.
Third, the $\bm{\nu}_i$ for $i \in \{n-n_c+1:n-1\}$ must be linearly independent.
The last two properties state that the non-trivial reactions
form simple neutral salts (i.e.~salts formed from two ions) and are independent.
These three properties imply that $\{\bm{\nu}_1, \ldots, \bm{\nu}_n, \bm{z} \}$
is a basis for $\R^n$ and therefore $\bm{Z}$ is invertible.

Through $\bm{Z}$ the electrochemical potentials
transform by
$\bm{\mu}_Z := \bm{Z} \bm{\mu}$.
Under electroneutrality, entries $(\bm{\mu}_Z)_i$ for $i \in \{1:n-1\}$ are independent of the
electrical state of the system and represent chemical potentials.
Indeed, for $i \in \{1:n-n_c\}$, $(\bm{\mu}_Z)_i = \mu_i$ is the
electrochemical potential of species $i$, but because species $i$ is uncharged, its electrochemical potential does
not depend on the electrical state and is just a chemical potential.
For $i \in \{n-n_c+1 : n - 1\}$, $(\bm{\mu}_Z)_i = \bm{\nu}_i^{\top} \bm{\mu}$
does not depend on the electrical state because it is the chemical potential of
the neutral salt formed in the reaction corresponding to $\bm{\nu}_i$
\cite{newman1965mass,van2023structural}.
Following \cite{van2023structural}, we decompose $\bm{\mu}_Z$ as
(recall that $F$ is the Faraday constant)
$\bm{\mu}_Z^{\top} = [\bm{\mu}_{\nu}^{\top} , (F \norm{\bm{z}} \Phi_Z)]$
where $\bm{\mu}_{\nu}$ is an $(n-1) \times 1$ column vector of the chemical
potentials and $\Phi_Z$ is a scalar-valued
\textit{salt-charge potential}. The potential $\Phi_Z$ can be related to the 
potential of a hypothetical probe reference electrode immersed in the mixture, and
serves to quantify the electrical state.
However, no material parameters in the model depend on $\Phi_Z$.

The concentrations and fluxes transform by
$\bm{c}_Z := \bm{Z}^{-\top} \bm{c}$ and $\bm{N}_Z := \bm{Z}^{-\top} \bm{N}$.
Since $\bm{\mu}_Z := \bm{Z} \bm{\mu}$ these transformations preserve
the structure of the volumetric Gibbs free energy
$\tilde{G} = \bm{c}^{\top} \bm{\mu} = \bm{c}_Z^{\top} \bm{\mu}_Z$
and energy dissipation function
$T \dot{s} = - \tr(\bm{N}^{\top} \nabla \bm{\mu}) 
= - \tr(\bm{N}_Z^{\top} \nabla \bm{\mu}_Z)$.
As in \cite{van2023structural} we decompose
\begin{equation} \label{eq:cZ_nZ_def}
	\bm{c}_Z = \begin{bmatrix}
		\bm{c}_{\nu} \\
		0
	\end{bmatrix}, \qquad
	\bm{N}_Z = \begin{bmatrix}
		\bm{N}_{\nu} \\
		J^{\top} / (F \norm{\bm{z}})
	\end{bmatrix}.
\end{equation}
The condition $(\bm{c}_Z)_n=0$ is equivalent to the electroneutrality condition
\cref{eq:v_ensosm_en}, and $J = F \bm{N}^{\top} \bm{z}$ is the current density.
The variables $\bm{c}_{\nu}$ can be thought of as the molar concentrations
of the salts formed in the neutralizing reactions, and 
$\bm{N}_{\nu}$ can be thought of as their molar fluxes.
Following \cite{van2023structural}, we refer to these hypothetical salts 
as \textit{components} as opposed to \textit{species}.
Using the equivalent electroneutrality condition $(\bm{c}_Z)_n=0$, one verifies 
that the governing equations \crefrange{eq:v_ensosm_mc}{eq:v_ensosm_osm}
transform as
\begin{alignat}{2}
	\frac{\partial}{\partial t}
	\begin{bmatrix}
		\bm{c}_{\nu} \\
		0
	\end{bmatrix}
	+ \begin{bmatrix}
		\nabla \cdot \bm{N}_{\nu} \\
		\nabla \cdot J
	\end{bmatrix} &= 0
	\quad &&\textrm{in } \Omega,
	\label{eq:vt_ensosm_mc} \\
	v &= \bm{\psi}_Z^{\top} \bm{N}_Z
	\quad &&\textrm{in } \Omega,
	\label{eq:vt_ensosm_mavg} \\
	-\nabla \bm{\mu}_Z + \bm{\psi}_Z \nabla p &= \bm{M}_Z \bm{N}_Z
	\quad &&\textrm{in } \Omega,
	\label{eq:vt_ensosm_osm}
\end{alignat}
where $\bm{\psi}_Z := \bm{Z} \bm{\psi}$ and 
$\bm{M}_Z := \bm{Z} \bm{M} \bm{Z}^{\top}$.
Since the \textit{transformed Onsager transport matrix} 
$\bm{M}_Z$ is congruent to $\bm{M}$,
it retains the crucial structural properties of
symmetry and positive-semidefiniteness (with a nullspace of dimension one). 

\Crefrange{eq:vt_ensosm_mc}{eq:vt_ensosm_osm} structurally resemble
the uncharged OSM equations.
In particular, the stationary version of 
\crefrange{eq:vt_ensosm_mc}{eq:vt_ensosm_osm}
(removing the time derivative term in \cref{eq:vt_ensosm_mc})
may be coupled to the stationary Stokes equations for $v$ and $p$,
and the resulting problem can be solved using the stationary SOSM
solver of \cite{baier2024high}.
The Picard linearized analysis of \cite{baier2024high} 
is applicable in such a setting.
However, instead of working with this formulation we will instead further
develop \cref{eq:vt_ensosm_osm} by expressing the chemical potential gradients in
terms of mole fraction gradients.
Consequently we will not need a constitutive formula for the chemical
potentials (or equivalently, activities) and will only need
a constitutive formula for the thermodynamic factors.
The mole fractions transform as
$\bm{x}_Z := \bm{Z}^{-\top} \bm{x}$ and electroneutrality implies
$(\bm{x}_Z)_n = 0$, hence we write
$\bm{x}_Z^{\top} = [\bm{x}_{\nu}^{\top} , 0]$
with $\bm{x}_{\nu}$ an $(n-1) \times 1$ column vector that parametrizes 
composition of the electroneutral mixture.
The normalization constraint $\sum_{i=1}^n x_i = 1$ transforms to
\begin{equation} \label{eq:en_mfs_normalization}
	\bm{\nu}_Z^{\top} \bm{x}_{\nu} = 1
\end{equation}
where $(\bm{\nu}_Z)_i = \sum_{j=1}^n \bm{Z}_{ij}$ for $i \in \{1 : n - 1\}$.
Generalizing \cite[eq.~(8.22)]{van2023structural} to the non-isobaric case,
we have
\begin{equation} \label{eq:vt_ensosm_gradmu}
	\nabla \bm{\mu}_{\nu} = 
	RT \bm{X}_{\nu}^0 \nabla \bm{x}_{\nu} + \bm{V}_{\nu} \nabla p
\end{equation}
where the $(n-1) \times (n-1)$ matrix $\bm{X}_{\nu}^0$ encodes
the thermodynamic factors under electroneutrality, and
$\bm{V}_{\nu}$ is an $(n-1) \times 1$ column vector of partial molar volumes:
\begin{equation*}
	\bm{V}_{\nu} := 
	\Big( 
	\frac{\partial \bm{\mu}_{\nu}}{\partial p}
	\Big)_{T, \bm{x}_{\nu}=\textrm{constant.}}
\end{equation*}
Partial molar volumes can be recovered from mass-density measurements
using Newman's formula
\cite[Appendix A]{newman2021electrochemical}.
Both $\bm{X}_{\nu}^0$ and $\bm{V}_{\nu}$ are functions of
$(T, p, \bm{x}_{\nu})$ only.

\subsection{Augmentation strategy}

As mentioned previously, the Onsager transport matrix $\bm{M}$ has a nullspace
of dimension one, and only $n-1$ of the OSM equations in \cref{eq:ensosm_osm}
are linearly independent.
It is the mass-average constraint \cref{eq:ensosm_mavg}, in conjunction
with the OSM equations \cref{eq:ensosm_osm}, that allows for the molar fluxes
to be uniquely determined.
To weakly enforce the mass-average constraint at the discrete level
and simultaneously address the fact that $\bm{M}$ is singular, we
employ the augmentation strategy introduced in
\cite{ern1994multicomponent,helfand1960inversion} and used in
\cite{aznaran2024finite,baier2024high,van2022augmented}.
In the present transformed setting,
we introduce a user-chosen augmentation parameter $\gamma > 0$
and add a multiple of the mass-average constraint \cref{eq:vt_ensosm_mavg} to
the OSM equation \cref{eq:vt_ensosm_osm} in the following way:
\begin{equation} \label{eq:vt_ensosm_aug_osm}
	-\nabla \bm{\mu}_Z + \bm{\psi}_Z \nabla p 
	+ \gamma \bm{\psi}_Z v
	= \gamma \bm{\psi}_Z \bm{\psi}_Z^{\top} \bm{N}_Z + \bm{M}_Z \bm{N}_Z
	= \bm{M}_Z^{\gamma} \bm{N}_Z
	\quad \textrm{in } \Omega,
\end{equation}
where $\bm{M}_Z^{\gamma} := \gamma \bm{\psi}_Z \bm{\psi}_Z^{\top} + \bm{M}_Z$
is a nonsingular augmented transport matrix.
Following \cite{aznaran2024finite,baier2024high}, for symmetry we also incorporate 
the mass-average constraint in \cref{eq:ensosm_ns} through%
\begin{equation} \label{eq:ensosm_aug_ns}
	\frac{\partial (\rho v)}{\partial t}
	+ \nabla \cdot (\rho v \otimes v)
	+ \nabla p
	+ \gamma (v - \bm{\psi}_Z^{\top} \bm{N}_Z)
	- \nabla \cdot \tau = \rho f
	\quad \textrm{in } \Omega.
\end{equation}
We also only explicitly enforce the divergence of the mass-average constraint,
i.e.~
\begin{equation} \label{eq:vt_ensosm_aug_mavg}
	\nabla \cdot v = \nabla \cdot (\bm{\psi}_Z^{\top} \bm{N}_Z)
	\quad \textrm{in } \Omega.
\end{equation}
We shall employ the augmented equations
\crefrange{eq:vt_ensosm_aug_osm}{eq:vt_ensosm_aug_mavg} because they
yield a system with the same number of equations as unknowns,
and, in an appropriate Picard linearized setting, yield a well-posed symmetric 
saddle point problem \cite{aznaran2024finite,baier2024high}.

\subsection{Full problem formulation} \label{sec:cts_final_eqns}
Our final formulation of the electroneutral NSOSM problem
comprises the following set of equations, taken from
\cref{eq:vt_ensosm_mc}, \cref{eq:ensosm_aug_ns}, \cref{eq:vt_ensosm_aug_mavg}
and combining \cref{eq:vt_ensosm_gradmu} with \cref{eq:vt_ensosm_aug_osm}.
The unknown functions of space and time to be solved for are
the velocity $v$,
transformed fluxes $\bm{N}_{\nu}$,
current density $J$,
pressure $p$,
transformed mole fractions $\bm{x}_{\nu}$
and salt-charge potential $\Phi_Z$, such that:%
\begin{subequations} \label{eq:final_ensosm}
\allowdisplaybreaks
\begin{alignat}{2}
	\frac{\partial (\rho v)}{\partial t}
	+ \nabla \cdot (\rho v \otimes v)
	+ \nabla p
	+ \gamma (v - \bm{\psi}_Z^{\top} \bm{N}_Z)
	- \nabla \cdot \tau &= \rho f
	\quad &&\textrm{in } \Omega, \label{eq:final_ns} \\
	-\begin{bmatrix}
		RT \bm{X}_{\nu}^0 &  0 \\
		0 & F \norm{\bm{z}}
	\end{bmatrix}
	\begin{bmatrix}
		\nabla \bm{x}_{\nu} \\
		\nabla \Phi_Z
	\end{bmatrix}
	+
	\bigg\{
	\bm{\psi}_Z - 
	\begin{bmatrix}
		\bm{V}_{\nu} \\
		0
	\end{bmatrix} \bigg\}
	\nabla p + \gamma \bm{\psi}_Z v
	&= \bm{M}_Z^{\gamma} \bm{N}_Z
	\quad &&\textrm{in } \Omega, \label{eq:final_osm} \\
	\nabla \cdot (v - \bm{\psi}_Z^{\top} \bm{N}_Z) &= 0
	\quad &&\textrm{in } \Omega,\label{eq:final_mavg} \\
	\frac{\partial}{\partial t}
	\begin{bmatrix}
		\bm{c}_{\nu} \\
		0
	\end{bmatrix}
	+ \begin{bmatrix}
		\nabla \cdot \bm{N}_{\nu} \\
		\nabla \cdot J
	\end{bmatrix} &= 0
	\quad &&\textrm{in } \Omega. \label{eq:final_mc}
\end{alignat}
\end{subequations}
This problem must be supplemented with suitable initial and boundary conditions,
which we discuss in forthcoming sections.

We recall how all quantities appearing in \cref{eq:final_ensosm} depend on the
unknowns.
The viscous stress $\tau$ in \cref{eq:final_ns} is given by
\cref{eq:ensosm_newtonian} and the viscosities in \cref{eq:ensosm_newtonian} are
assumed to be known functions of $(T, p, \bm{x}_{\nu})$.
The vector $\bm{N}_Z$ was defined in \cref{eq:cZ_nZ_def} and is simply a
concatenation of the fluxes $\bm{N}_{\nu}$ and scaled current density $J$.
The constant $\gamma > 0$ is a user-chosen augmentation parameter
and $f$ is a known forcing term.
All remaining quantities in \cref{eq:final_ensosm} are assumed to be known
Lipschitz continuous functions of $(T, p, \bm{x}_{\nu})$.
The density $\rho$, concentrations $\bm{c}_{\nu}$ and partial molar volumes
$\bm{V}_{\nu}$ can be determined as functions of $(T, p, \bm{x}_{\nu})$
using the volumetric EOS.
The vector $\bm{\psi}_Z = \bm{Z} \bm{m} / \rho$ depends on $\rho$ only.
The matrix $\bm{X}_{\nu}^0$ of thermodynamic factors can be expressed
in terms of $(n-1)(n-2)/2$ independent Darken factors,
while the transport matrix $\bm{M}_Z^{\gamma}$ can be expressed in
terms of $n(n-1)/2$ independent Stefan--Maxwell diffusivities
\cite{van2023structural}. The dependence of these material parameters on
$(T, p, \bm{x}_{\nu})$ is typically modelled by fitting experimental data.

\section{Spatial discretization} \label{sec:spatial_disc}

In this section, we introduce finite element methods to spatially discretize the
electroneutral NSOSM equations \cref{eq:final_ensosm}
in steady form:
\begin{subequations} \label{eq:steady_ensosm}
\allowdisplaybreaks
\begin{alignat}{2}
	\nabla \cdot (\rho v \otimes v)
	+ \nabla p
	+ \gamma (v - \bm{\psi}_Z^{\top} \bm{N}_Z)
	- \nabla \cdot \tau &= \rho f
	\quad &&\textrm{in } \Omega, \label{eq:steady_ns} \\
	-\begin{bmatrix}
		RT \bm{X}_{\nu}^0 &  0 \\
		0 & F \norm{\bm{z}}
	\end{bmatrix}
	\begin{bmatrix}
		\nabla \bm{x}_{\nu} \\
		\nabla \Phi_Z
	\end{bmatrix}
	+
	\bigg\{
	\bm{\psi}_Z - 
	\begin{bmatrix}
		\bm{V}_{\nu} \\
		0
	\end{bmatrix} \bigg\}
	\nabla p + \gamma \bm{\psi}_Z v
	&= \bm{M}_Z^{\gamma} \bm{N}_Z
	\quad &&\textrm{in } \Omega, \label{eq:steady_osm} \\
	\nabla \cdot (v - \bm{\psi}_Z^{\top} \bm{N}_Z) &= 0
	\quad &&\textrm{in } \Omega,\label{eq:steady_mavg} \\
	\nabla \cdot \bm{N}_{\nu} &= 0
	\quad &&\textrm{in } \Omega, \label{eq:steady_mc} \\
	\nabla \cdot J &= 0
	\quad &&\textrm{in } \Omega. \label{eq:steady_kl}
\end{alignat}
\end{subequations}
The steady problem \cref{eq:steady_ensosm} must be supplemented with suitable 
boundary conditions and integral constraints, which we now describe.

\subsection{Boundary conditions}
Let $\Gamma := \partial \Omega$ denote the boundary of $\Omega$
and $n_{\Gamma}$ the unit outward normal on $\Gamma$.
We consider the following flux\footnote{We do not consider Dirichlet BCs 
in this section, but see \cref{sec:rotating_electrode} for an example
of how they can be implemented.}
boundary conditions:
\begin{subequations} \label{eq:steady_ensosm_bcs}
\allowdisplaybreaks
\begin{alignat}{2}
	v &= 
	[(\bm{\psi}_Z^{\top} \bm{N}_Z) \cdot n_{\Gamma}] n_{\Gamma}
	+ g_{v_{\parallel}}
	&&\quad \textrm{on } \Gamma, \label{eq:steady_ensosm_bcs_v} \\
	(\bm{N}_{\nu})_i \cdot n_{\Gamma} &= g_i \label{eq:steady_ensosm_bcs_N}
	&&\quad \textrm{on } \Gamma \quad \forall i \in \{1 : n-1\}, \\
	J \cdot n_{\Gamma} &= g_J &&\quad \textrm{on } \Gamma.
	\label{eq:steady_ensosm_bcs_J}
\end{alignat}
\end{subequations}
The functions $g_{v_{\parallel}} : \Gamma \to \R^{d},\ g_i : \Gamma \to \R$
and $g_J : \Gamma \to \R$ are prescribed data.
We assume $g_{v_{\parallel}} \cdot n_{\Gamma} = 0$ on $\Gamma$, so that
\cref{eq:steady_ensosm_bcs_v} enforces the mass-average constraint
\cref{eq:vt_ensosm_mavg} in the normal direction and
$v = g_{v_{\parallel}}$ in the tangential directions.

In practical applications the prescribed normal fluxes $g_i$ and $g_J$
may be known algebraic functions of the state variables
$(T, p, \bm{x}_{\nu})$ and $\Phi_Z$.
We shall allow for such dependencies.
An example is that of a Butler--Volmer boundary condition
\cite{butler1924studies,dickinson2020butler,erdey1930theorie,
newman2021electrochemical}, which on an electrode interface
$\Gamma_{\textrm{e}} \subset \Gamma$ relates the normal
component of current density to the \textit{overpotential}.
The overpotential quantifies the electrical potential difference
across the interface.
Butler--Volmer BCs generally depend on $(T, p, \bm{x}_{\nu})$ and $\Phi_Z$ 
nonlinearly, and we consider this case in \cref{sec:numerical}.
However, if the overpotential is sufficiently small and its dependence on 
electrolyte composition is negligible, it is worth noting that
Butler--Volmer BCs may be approximated by
a linear relation \cite{newman2021electrochemical}
\begin{equation} \label{eq:bv_linear}
	g_J = -i_0
	\frac{n_e F}{RT} 
	(V_{\textrm{e}} - \Phi_Z)
	\quad \textrm{on } \Gamma_{\textrm{e}}.
\end{equation}
Here $i_0$ denotes the \textit{exchange-current density},
which generally depends on $(T, p, \bm{x}_{\nu})$. Furthermore,
$V_e$ denotes the electrode voltage,
and $n_e$ is a stoichiometric coefficient representing the number of electrons 
transferred in the electrode reaction.

\subsection{Integral constraints} \label{eq:steady_constraints}

In general, the steady problem \cref{eq:steady_ensosm} with
boundary conditions \cref{eq:steady_ensosm_bcs} is not uniquely solvable, and 
$1 \leq k \leq n+1$ additional constraints
must be imposed for uniqueness.
The necessity of such constraints in the uncharged, steady SOSM setting
was pointed out in \cite{baier2024high}.
However, the present situation is more complicated than that of 
\cite{baier2024high} because we allow for the boundary data in 
\crefrange{eq:steady_ensosm_bcs_N}{eq:steady_ensosm_bcs_J}
to depend on the local state of the mixture.
Further complications will arise in \cref{sec:temporal_disc} when
we consider the transient problem, but for now we focus on the steady case.

The need for constraints can be motivated as follows.
Integrate the solenoidal condition \cref{eq:steady_kl} on $J$
over $\Omega$ and use the divergence theorem together with the BC
\cref{eq:steady_ensosm_bcs_J} to obtain a fundamental
\textit{compatibility condition} on the problem data:
\begin{equation} \label{eq:compat_cond_J}
	\int_{\Gamma} g_J \mathop{\mathrm{d}\Gamma} = 0,
\end{equation}
which states that the total charge flux into the electroneutral mixture is zero.
By contrast, in a general electrolytic cell the total component molar fluxes need 
not be zero, unless the mixture is at a steady state.
In the steady case, however, the continuity equations
\cref{eq:steady_mc} and BCs
\cref{eq:steady_ensosm_bcs_N} yield $n-1$ additional compatibility conditions:
\begin{equation} \label{eq:compat_conds_N}
	\int_{\Gamma} g_i \mathop{\mathrm{d}\Gamma} = 0
	\quad \forall i \in \{1 : n - 1\}.
\end{equation}

For the steady problem \cref{eq:steady_ensosm} with BCs \cref{eq:steady_ensosm_bcs}
to admit a solution, the conditions in \cref{eq:compat_cond_J,eq:compat_conds_N} 
must be satisfiable.
Assuming this is so, let $l \leq n$ denote the number of independent constraints
that \cref{eq:compat_cond_J,eq:compat_conds_N} impose on the unknowns
$(p, \bm{x}_{\nu}, \Phi_Z)$ due to their (possible) appearance in $g_J$ and $g_i$.
Finally, similarly integrate \cref{eq:steady_mavg} over $\Omega$ and
use the BC \cref{eq:steady_ensosm_bcs_v} to obtain an $(n+1)$-th
compatibility condition
\begin{equation} \label{eq:compat_conds_v}
	\int_{\Gamma} 
	g_{v_{\parallel}} \cdot n_{\Gamma} \mathop{\mathrm{d}\Gamma} = 0,
\end{equation}
which holds automatically since $g_{v_{\parallel}} \cdot n_{\Gamma} = 0$.
Hence, we have shown that integrating the $n+1$ equations in
\crefrange{eq:steady_mavg}{eq:steady_kl} over $\Omega$ leads to only $l \leq n$
independent constraints on the solution.
An additional $k = n + 1 - l$ constraints must therefore be imposed for uniqueness.
The same argument will apply at the discrete level.
Namely, if $D$ denotes the number of discrete degrees of freedom, then
discretization of \crefrange{eq:steady_ensosm}{eq:steady_ensosm_bcs} will result
in a system of $D$ equations, but only $D-k$ of these will be independent.
This motivates why $k$ additional constraints are needed for uniqueness.
We will give examples of choices of constraints in \cref{sec:constraints_discrete}.
Note that these constraints do not need to be integral 
constraints, but they usually are in practice.

\subsection{Finite element spaces} \label{sec:fem_spaces}

We use the standard notation for Lebesgue and Sobolev spaces
$L^2(\Omega),\ W^{1, \infty}(\Omega),\ H^1(\Omega),\ H(\div; \Omega)$
and their norms \cite{ern2021finiteI}.
We use $(\cdot, \cdot)_{\Omega}$ and $\langle\cdot, \cdot\rangle_{\Gamma}$
to denote the $L^2$-inner products on $\Omega$ and $\Gamma$ respectively.
Moreover, we write
$L_0^2(\Omega) := \{ q \in L^2(\Omega) : (q, 1)_{\Omega} = 0 \}$
for the subspace of functions in $L^2(\Omega)$ with vanishing mean.
We write
$H_0^1(\Omega) := \{ u \in H^1(\Omega) : u|_{\Gamma} = 0 \}$ and
$H_0(\div; \Omega) := \{ u \in H(\div; \Omega) : 
(u \cdot n_{\Gamma})|_{\Gamma} = 0 \}$ for subspaces with vanishing traces.

To discretize \cref{eq:steady_ensosm} we introduce finite dimensional
finite element subspaces that may depend on a parameter $h \in (0, \infty)$
representing, for example, the mesh size:
\begin{subequations} \label{eq:fem_spaces}
\allowdisplaybreaks
\begin{alignat}{2}
	&V_h \subset H^1(\Omega)^d, &&\qquad V_{0h} := V_h \cap H_0^1(\Omega)^d, \\
	&P_h \subset L^2(\Omega), &&\qquad P_{0h} := P_h \cap L_0^2(\Omega), \\
	&N_h \subset H(\div;\Omega), &&\qquad N_{0h} := N_h \cap H_0(\div;\Omega), \\
	&X_h \subset L^2(\Omega), &&\qquad X_{0h} := X_h \cap L_0^2(\Omega).
\end{alignat}
\end{subequations}
Importantly, we assume $P_h$ and $X_h$ contain the constant functions,
i.e.~$1 \in P_h \cap X_h$.

We shall seek the discrete velocity $v_h \in V_h$ and pressure $p_h \in P_h$.
We assume that $(V_h, P_h)$ forms an inf-sup stable Stokes pair
\cite{ern2021finiteII}, in the sense that:
\begin{equation} \label{eq:fem_spaces_a1}
	\beta \norm{q_h}_{L^2(\Omega)} \leq
	\sup_{u_h \in V_{0h}} 
	\frac{(q_h, \div u_h)_{\Omega}}{\norm{u_h}_{H^1(\Omega)^d}}
	\qquad
	\forall q_h \in P_{0h},
\end{equation}
for some $\beta > 0$ independent of $h$.
We shall seek the discrete fluxes $(\bm{N}_{\nu,h})_i \in N_h$ for
$i \in \{1 : n - 1\}$ and current density $J_h \in N_h$.
We seek the discrete mole fractions
$(\bm{x}_{\nu,h})_i \in X_h$ for $i \in \{1 : n - 1\}$
and salt-charge potential $\Phi_{Z,h} \in X_h$.
We assume that $(N_h, X_h)$ forms a divergence-free and inf-sup stable
mixed-Poisson pair \cite{ern2021finiteII}, i.e.~that:
\begin{equation} \label{eq:fem_spaces_a2}
	\div N_{0h} \subset X_{0h} \quad \textrm{and} \quad
	\beta' \norm{y_h}_{L^2(\Omega)} \leq
	\sup_{K_h \in N_{0h}}
	\frac{(y_h, \div K_h)_{\Omega}}{\norm{K_h}_{H(\div;\Omega)}}
	\qquad
	\forall y_h \in X_{0h},
\end{equation}
for some $\beta' > 0$ independent of $h$.
Assumptions \cref{eq:fem_spaces_a1,eq:fem_spaces_a2} are motivated by the
analysis of \cite{baier2024high}, where similar assumptions are shown
to ensure well-posedness of a Picard linearized SOSM system.

Following \cite{baier2024high}, in our numerical experiments
we employ the degree $k \geq 2$ Taylor--Hood pair
\cite{boffi1997three,taylor1973numerical} for $(V_h, P_h)$
(but cf.\ \cite{baier2024high} for a discussion of other choices,
such as Scott--Vogelius \cite{scott1985norm}).
For $(N_h, X_h)$ we employ either
the $\mathbb{BDM}_k$--$\mathbb{DG}_{k-1}$ pair
\cite{brezzi1985two,nedelec1986new}
or the $\mathbb{RT}_k$--$\mathbb{DG}_{k-1}$
\cite{raviart1977mixed} pair.
These spaces are standard in finite element software packages,
extend to high orders, and are applicable in two or three
spatial dimensions on triangular, tetrahedral, quadrilateral or hexahedral 
(possibly curved) meshes.

\subsection{Lipschitz continuous reconstruction operators}

Our discretization in \cref{sec:fem_disc}
will involve integration-by-parts of the gradient
terms on the left side of \cref{eq:steady_osm}.
This requires the entries of $\bm{X}_{\nu}^0$, $\bm{\psi}_{Z}$ and
$\bm{V}_{\nu}$ to lie in $W^{1, \infty}(\Omega)$.
Discretization of \cref{eq:steady_mavg} will likewise require
$\bm{\psi}_{Z} \in (W^{1, \infty}(\Omega))^n$.
Recall that we assume these quantities to be known Lipschitz continuous
functions of $(T, p, \bm{x}_{\nu})$.
However, at the discrete level,
the spaces $P_h$ or $X_h$ may be discontinuous.
For this reason we shall evaluate thermodynamic properties such as 
$\bm{X}_{\nu}^0$, $\bm{\psi}_{Z}$, $\bm{V}_{\nu}$
using Lipschitz continuous reconstructions of $p_h$ and $\bm{x}_{\nu,h}$.
To be precise, we assume that smoothing operators
\begin{equation*}
	\pi_{P_h} : P_h \to P_h \cap W^{1, \infty}(\Omega)
	\quad \textrm{ and } \quad
	\pi_{X_h} : X_h \to X_h \cap W^{1, \infty}(\Omega)
\end{equation*}
are available. In our numerical simulations we take
$\pi_{P_h}$ to be the $L^2(\Omega)$-projection of $P_h$ into 
$P_h \cap W^{1, \infty}(\Omega)$ and likewise for $X_h$.
Nodal averaging operators could alternatively be employed
\cite{ern2017finite}.
Given $p_h$ and $\bm{x}_{\nu,h}$ we introduce their \textit{reconstructions}
\begin{equation*}
	\widetilde{p_h} := \pi_{P_h} p_h
	\quad \textrm{ and } \quad
	\widetilde{(\bm{x}_{\nu,h})}_i 
	:= \pi_{X_h} (\bm{x}_{\nu,h})_i
	\Big/
	\textstyle\sum_{j=1}^{n-1} (\bm{\nu}_Z)_j \cdot \pi_{X_h}(\bm{x}_{\nu,h})_j.
\end{equation*}
Note that the reconstructed mole fractions are normalized to satisfy condition 
\cref{eq:en_mfs_normalization} exactly.
We then write
$\widetilde{\bm{X}_{\nu}^0}$, $\widetilde{\bm{\psi}_{Z}}$,
$\widetilde{\bm{V}_{\nu}}$, $\widetilde{\rho}$, and so on,
to denote quantities
evaluated with the reconstructed $\widetilde{p_h}$ and
$\widetilde{\bm{x}_{\nu,h}}$ instead of $p_h$ and $\bm{x}_{\nu,h}$.

\subsection{Discretized problem} \label{sec:fem_disc}

Our discrete variational formulation of \cref{eq:steady_ensosm} can be obtained
by multiplying the equations with suitable test functions and integrating
over $\Omega$.
The pressure and viscous terms in \cref{eq:steady_ns}, as well as
the gradient terms on the left side of \cref{eq:steady_osm}, are 
integrated by parts, and all boundary terms vanish owing to our BCs
in \cref{eq:steady_ensosm_bcs}.
Following \cite{baier2024high} we also add \textit{density consistency terms} to 
\cref{eq:steady_mavg}.

The precise discretization we consider is as follows.
We seek discrete functions $v_h \in V_h$, $p_h \in P_h$,
$\bm{N}_{\nu,h} \in (N_h)^{n-1}$,
$J_h \in N_h$,
$\bm{x}_{\nu,h} \in (X_h)^{n-1}$
and
$\Phi_{Z,h} \in X_h$.
Let
\begin{equation*}
	\bm{N}_{Z,h} := 
	\begin{bmatrix}
		\bm{N}_{\nu,h} \\
		J_h^{\top} / (F \norm{z})
	\end{bmatrix} \in (N_h)^n,
\end{equation*}
which is analogous to $\bm{N}_Z$ in \cref{eq:cZ_nZ_def}.
The discrete variational problem reads:
\begin{subequations} \label{eq:disc_steady_ensosm}
\allowdisplaybreaks
\begin{align}
\begin{split}
	&\big( \nabla \cdot (\widetilde{\rho} v_h \otimes v_h)
	+ \gamma (v_h - \widetilde{\bm{\psi}_Z^{\top}} \bm{N}_{Z,h})
	- \widetilde{\rho} f, u_h \big)_{\Omega}
	- \big( p_h, \nabla \cdot u_h \big)_{\Omega} \\
	&+ \ \big( 2 \widetilde{\eta} \epsilon(v_h) + 
	(\widetilde{\zeta} - 2 \widetilde{\eta} / d) (\nabla \cdot v_h) \Id,  
	\epsilon(u_h) \big)_{\Omega}
	= 0
	\quad \forall u_h \in V_{0h},
\end{split}
\\
\begin{split}
	&\Bigg(
	\begin{bmatrix}
		\bm{x}_{\nu,h} \\
		\Phi_{Z,h}
	\end{bmatrix},
	\nabla \cdot
	\begin{bmatrix}
		RT \widetilde{\bm{X}_{\nu}^0} \bm{W}_h \\
		F \norm{\bm{z}} K_h^{\top}
	\end{bmatrix}
	\Bigg)_{\Omega}
	-
	\Bigg(p_h,
	\nabla \cdot
	\Bigg\{
	\widetilde{\bm{\psi}_{Z}^{\top}}
	\begin{bmatrix}
		\bm{W}_h \\
		K_h^{\top}
	\end{bmatrix}
	- \widetilde{\bm{V}_{\nu}^{\top}} \bm{W}_h
	\Bigg\}
	\Bigg)_{\Omega} \\
	&+ \ \Bigg(
	\gamma  \widetilde{\bm{\psi}_{Z}} v_h,
	\begin{bmatrix}
		\bm{W}_h \\
		K_h^{\top}
	\end{bmatrix}
	\Bigg)_{\Omega}
	=
	\Bigg(
	\widetilde{\bm{M}_Z^{\gamma}} \bm{N}_{Z,h},
	\begin{bmatrix}
		\bm{W}_h \\
		K_h^{\top}
	\end{bmatrix} \Bigg)_{\Omega}
	\quad \forall
	\begin{bmatrix}
		\bm{W}_h \\
		K_h^{\top}
	\end{bmatrix} \in (N_{0h})^n,
\end{split} \label{eq:disc_steady_ensosm_osm} \\
\begin{split}
	\Big( \nabla \cdot (v_h - \widetilde{\bm{\psi}_Z^{\top}} 
	\bm{N}_{Z,h}), q_h \Big)_{\Omega}
	- \Big\langle
	n_{\Gamma} \cdot (v_h - \widetilde{\bm{\psi}_Z^{\top}} \bm{N}_{Z,h}), q_h
	\Big\rangle_{\Gamma} = 0
	\quad \forall q_h \in P_h,
\end{split} \label{eq:disc_steady_ensosm_mavg}  \\
\begin{split}
	\big(
	\nabla \cdot
	\bm{N}_{Z,h},
	\bm{y}_h
	\big)_{\Omega} = 0
	\quad \forall \bm{y}_h \in (X_h)^n.
\end{split} \label{eq:disc_steady_ensosm_cty}
\end{align}
\end{subequations}
Moreover, we strongly impose the following discrete analogue
of the BCs in \cref{eq:steady_ensosm_bcs}:
\begin{subequations} \label{eq:disc_steady_bcs}
\allowdisplaybreaks
\begin{alignat}{2}
	v_h &= \pi_{V_h} \Big(
	[(\widetilde{\bm{\psi}_Z^{\top}} \bm{N}_{Z,h}) \cdot n_{\Gamma}] n_{\Gamma}
	+ g_{v_{\parallel}} \Big)
	&&\quad \textrm{on } \Gamma,
	\label{eq:disc_steady_bcs_v} \\
	(\bm{N}_{\nu,h})_i \cdot n_{\Gamma} &= \pi_{N_h} \big( \widetilde{g_i} \big)
	&&\quad \textrm{on } \Gamma \quad \forall i \in \{1 : n-1\},
	\label{eq:disc_steady_bcs_N} \\
	J_h \cdot n_{\Gamma} &= 
	\pi_{N_h} \big(
	\widetilde{g_J} \big) &&\quad \textrm{on } \Gamma.
	\label{eq:disc_steady_bcs_J}
\end{alignat}
\end{subequations}
Here $\pi_{V_h}$ and $\pi_{N_h}$ are $L^2(\Gamma)$-projection operators
onto the discrete trace spaces%
\footnote{
For \cref{eq:proj_trace_N} to be well-defined, we implicitly assume 
$(K_h \cdot n_{\Gamma})|_{\Gamma} \in L^2(\Gamma) \ \forall K_h \in N_h$.
In practice, the finite element space $N_h$ will consist of piecewise smooth
functions, so that this indeed holds.
}
\begin{subequations}
\allowdisplaybreaks
\begin{align}
	\pi_{V_h} &: 
	[L^2(\Gamma)]^d \to \{ u_h|_{\Gamma} : u_h \in V_h \}
	\subset [L^2(\Gamma)]^d, \\
	\pi_{N_h} &: 
	L^2(\Gamma) \to \{ (K_h \cdot n_{\Gamma})|_{\Gamma} : K_h \in N_h \} \subset 
	L^2(\Gamma). \label{eq:proj_trace_N}
\end{align}
\end{subequations}
Conditions \crefrange{eq:disc_steady_ensosm}{eq:disc_steady_bcs} define our
discrete scheme.

\subsection{Integral constraints in the discrete setting}
\label{sec:constraints_discrete}

As already discussed in \cref{eq:steady_constraints}, for uniqueness of a solution 
to the discretized problem \crefrange{eq:disc_steady_ensosm}{eq:disc_steady_bcs},
additional constraints must be supplied.
To see this at the discrete level, let
$D := \dim \big( V_h \times P_h \times (N_h)^n \times (X_h)^n \big)$
and note that, once a finite element basis has been chosen,
the problem \crefrange{eq:disc_steady_ensosm}{eq:disc_steady_bcs} amounts to a
nonlinear system of $D$ equations in $D$ unknowns.
However, notice that \cref{eq:disc_steady_ensosm_mavg} holds automatically
when $q_h$ is a constant, as verified from integration by parts
(analogous to the way \cref{eq:compat_conds_v} holds automatically).
This property is one way of motivating density consistency terms
(i.e.~the boundary terms) in \cref{eq:disc_steady_ensosm_mavg} 
\cite{baier2024high}.
Also, taking the entries of $\bm{y}_h$ to be constants in
\cref{eq:disc_steady_ensosm_cty}, integrating by parts and using
\crefrange{eq:disc_steady_bcs_N}{eq:disc_steady_bcs_J} yields
(in analogy to \cref{eq:compat_cond_J,eq:compat_conds_N})%
\begin{equation} \label{eq:compat_conds_disc}
	\int_{\Gamma} \pi_{N_h} \big( \widetilde{g_i} \big)
	\mathop{\mathrm{d}\Gamma} = 0
	\quad \forall i \in \{1 : n - 1\} \quad \textrm{and} \quad
	\int_{\Gamma} \pi_{N_h} \big(
	\widetilde{g_J} \big) \mathop{\mathrm{d}\Gamma} = 0.
\end{equation}
Let $l \leq n$ denote the number of independent constraints that 
\cref{eq:compat_conds_disc} imposes on the discrete solution.
Problem \crefrange{eq:disc_steady_ensosm}{eq:disc_steady_bcs} then consists of
only $D - k$ independent equations, where $k = n + 1 - l$.
Hence, an additional $k$ constraints must be imposed.

Since $k \geq 1$, at least one constraint is always required.
This constraint does not encode physical information, and instead
reflects our choice to solve for $n-1$ transformed mole fractions
despite the fact that only $n-2$ of them are independent.
Following \cite{baier2024high} we impose the
normalization condition \cref{eq:en_mfs_normalization} on average over $\Omega$:
\begin{equation} \label{eq:mf_constraint}
	\int_{\Omega} (\bm{\nu}_Z^{\top} \bm{x}_{\nu,h} - 1) 
	\mathop{\mathrm{d}\Omega}=0.
\end{equation}
Together with \cref{eq:mf_constraint}, the remaining $k-1$ constraints must be 
chosen on a case-by-case basis.
What choices of constraints are appropriate will depend on BCs and the functional 
dependence of the thermodynamic
properties (i.e.~the properties in \cref{eq:disc_steady_ensosm} with a tilde)
on $p$ and $\bm{x}_{\nu}$.
We now give examples of this for concreteness.

\textit{Case (i): The BCs in \cref{eq:steady_ensosm_bcs_N,eq:steady_ensosm_bcs_J}
do not depend on $p$, $\bm{x}_{\nu}$ or $\Phi_Z$}.
In this case, the compatibility conditions \cref{eq:compat_conds_disc} do not
impose any constraints on the solution, so $l=0$ and $n$ additional constraints
are required.
Since no thermodynamic properties depend on $\Phi_Z$,
this field is then undetermined up to an additive constant,
which can be fixed through a constraint such as
$\int_{\Omega} \Phi_{Z,h} \mathop{\mathrm{d}\Omega}=0$.
If no thermodynamic properties depend on $p$ then it too is
undetermined up to an additive constant, which can be fixed by means of
$\int_{\Omega} p_h \mathop{\mathrm{d}\Omega}=0$.
However, if any of thermodynamic properties depend on $p$ then additive
shifts of $p_h$ by a constant affect the physical content of the model.
In this case, constraints of the form
$\int_{\Omega} (p_h - \bar{p}) \mathop{\mathrm{d}\Omega}=0$
may be used, where $\bar{p} \in \R$ is a user-prescribed mean pressure.
If the volumetric EOS depends on $p$
(equivalently if $\rho$ depends on $p$) then an alternative constraint that is more
harmonious with experimentally available information may be
$\int_{\Omega} \widetilde{\rho} \mathop{\mathrm{d}\Omega} = M^{\textrm{tot}}$,
where $M^{\textrm{tot}} \in \R$ is the user-prescribed total mass of fluid in 
$\Omega$.
The final $n-2$ constraints may be chosen to express the relative abundances of 
the components in $\Omega$. For example, one may use
$\int_{\Omega} \widetilde{(\bm{c}_{\nu})}_i \mathop{\mathrm{d}\Omega} = 
c_i^{\textrm{tot}}$
for $i \in S$, where $S \subset \{1 : n-1\}$ contains $n-2$ indices, and 
$c_i^{\textrm{tot}} \in \R$ are user-prescribed total numbers of moles 
of the components in $\Omega$.

\textit{Case (ii): The BCs in \cref{eq:steady_ensosm_bcs_N,eq:steady_ensosm_bcs_J}
depend on the local state of the mixture}.
To illustrate how solution-dependent BCs fit into this framework,
consider the case of an electrolytic cell where $g_J$ follows a 
Butler--Volmer-type relation in \cref{eq:bv_linear} or some nonlinear analogue of 
this, while $g_i = \alpha_i g_J / F$ on $\Gamma$ for $i \in \{1 : n - 1\}$ with
$\{\alpha_i\}_{i=1}^{n-1}$ a set of known constants.
These BCs model electrode kinetics, and the $\alpha_i$ relate to the 
stoichiometric coefficients of the electrode reaction%
\footnote{%
Generally the value of $\alpha_i$ can differ on the anode and cathode; we assume 
here it is the same on both (this holds, for example, in symmetric cells and 
lithium-ion battery cells).
The relation $g_i = \alpha_i g_J / F$ is also valid on the insulating walls of the 
cell since therein both $g_i$ and $g_J$ are zero.
}.
Since each $g_i$ is a multiple of $g_J$, the compatibility
conditions in \cref{eq:compat_conds_disc} impose exactly one constraint
on the solution, so $l=1$ and $n-1$ additional constraints are required.
Since $\Phi_Z$ appears in the BC \cref{eq:bv_linear} it is no longer necessary (or
physically reasonable) to fix
$\int_{\Omega} \Phi_{Z,h} \mathop{\mathrm{d}\Omega}=0$.
Instead, the $n-1$ constraints should be chosen to determine the pressure
and relative component abundances, as in case (i).

\section{Temporal discretization} \label{sec:temporal_disc}

Our spatial discretization from 
\cref{sec:spatial_disc} can be applied in the transient setting using the method 
of lines.
However, special care must be taken to address the subtle interplay between
the component mass continuity equations, boundary conditions, volumetric EOS and 
integral constraints.
\subsection{Semi-discrete problem} \label{sec:semi_discrete}

Discretization in space (but not time) of the transient problem 
\cref{eq:final_ensosm} yields the following semi-discrete analogue of 
\cref{eq:disc_steady_ensosm}.
We seek time-dependent discrete functions $v_h \in V_h$, $p_h \in P_h$,
$\bm{N}_{\nu,h} \in (N_h)^{n-1}$,
$J_h \in N_h$,
$\bm{x}_{\nu,h} \in (X_h)^{n-1}$
and
$\Phi_{Z,h} \in X_h$ such that:
\begin{subequations} \label{eq:disc_transient_ensosm}
\allowdisplaybreaks
	\begin{align}
		\begin{split}
			&\frac{\mathrm{d}}{\mathrm{d} t}
			\big( \widetilde{\rho} v_h, u_h \big)_{\Omega}
			+ \ \big( \nabla \cdot (\widetilde{\rho} v_h \otimes v_h)
			+ \gamma (v_h - \widetilde{\bm{\psi}_Z^{\top}} \bm{N}_{Z,h})
			- \widetilde{\rho} f, u_h \big)_{\Omega} \\
			&- \ \big( p_h, \nabla \cdot u_h \big)_{\Omega}
			+ \ \big( 2 \widetilde{\eta} \epsilon(v_h) + 
			(\widetilde{\zeta} - 2 \widetilde{\eta} / d) (\nabla \cdot v_h) \Id,  
			\epsilon(u_h) \big)_{\Omega}
			= 0
			\quad \forall u_h \in V_{0h},
		\end{split}
		\\
		\begin{split}
			&\Bigg(
			\begin{bmatrix}
				\bm{x}_{\nu,h} \\
				\Phi_{Z,h}
			\end{bmatrix},
			\nabla \cdot
			\begin{bmatrix}
				RT \widetilde{\bm{X}_{\nu}^0} \bm{W}_h \\
				F \norm{\bm{z}} K_h^{\top}
			\end{bmatrix}
			\Bigg)_{\Omega}
			-
			\Bigg(p_h,
			\nabla \cdot
			\Bigg\{
			\widetilde{\bm{\psi}_{Z}^{\top}}
			\begin{bmatrix}
				\bm{W}_h \\
				K_h^{\top}
			\end{bmatrix}
			- \widetilde{\bm{V}_{\nu}^{\top}} \bm{W}_h
			\Bigg\}
			\Bigg)_{\Omega} \\
			&+ \ \Bigg(
			\gamma  \widetilde{\bm{\psi}_{Z}} v_h,
			\begin{bmatrix}
				\bm{W}_h \\
				K_h^{\top}
			\end{bmatrix}
			\Bigg)_{\Omega}
			=
			\Bigg(
			\widetilde{\bm{M}_Z^{\gamma}} \bm{N}_{Z,h},
			\begin{bmatrix}
				\bm{W}_h \\
				K_h^{\top}
			\end{bmatrix} \Bigg)_{\Omega}
			\quad \forall
			\begin{bmatrix}
				\bm{W}_h \\
				K_h^{\top}
			\end{bmatrix} \in (N_{0h})^n,
		\end{split} \\
		\begin{split}
			\Big( \nabla \cdot (v_h - \widetilde{\bm{\psi}_Z^{\top}} 
			\bm{N}_{Z,h}), q_h \Big)_{\Omega}
			- \Big\langle
			n_{\Gamma} \cdot (v_h - \widetilde{\bm{\psi}_Z^{\top}} \bm{N}_{Z,h}), 
			q_h
			\Big\rangle_{\Gamma} = 0
			\quad \forall q_h \in P_h,
		\end{split} \label{eq:disc_transient_ensosm_mavg}  \\
		\begin{split}
			\frac{\mathrm{d}}{\mathrm{d} t}
			\Bigg(
			\begin{bmatrix}
				\widetilde{\bm{c}_{\nu,h}} \\
				0
			\end{bmatrix},
			\bm{y}_h
			\Bigg)_{\Omega}
			+
			\big(
			\nabla \cdot
			\bm{N}_{Z,h},
			\bm{y}_h
			\big)_{\Omega} = 0
			\quad \forall \bm{y}_h \in (X_h)^n.
		\end{split} \label{eq:disc_transient_ensosm_cty}
	\end{align}
\end{subequations}
We supplement this problem with the same strongly enforced BCs 
\cref{eq:disc_steady_bcs} from the steady case,
and we permit the boundary data to depend on time.

\subsection{Integral constraints in the transient setting}
\label{sec:transient_constraints}

As in the steady case, problem \cref{eq:disc_transient_ensosm} with
BCs \cref{eq:disc_steady_bcs} requires integral constraints for well-posedness.
However, the interplay between the BCs, volumetric EOS and mass continuity 
equations becomes especially important in the transient setting, as we now
elucidate through considerations similar to those in 
\cref{sec:constraints_discrete}.

First, note that \cref{eq:disc_transient_ensosm_mavg} holds automatically when 
$q_h = 1$, as in the steady case. This suggests that at least one
integral constraint will be needed for well-posedness.
Next, taking constant entries of $\bm{y}_h$ in
\cref{eq:disc_transient_ensosm_cty} yields, similarly to 
\cref{eq:compat_conds_disc}, that
\begin{subequations} \label{eq:transient_compat_conds_disc}
\allowdisplaybreaks
\begin{align}
	\frac{\mathrm{d}}{\mathrm{d} t}
	\int_{\Omega} (\widetilde{\bm{c}_{\nu,h}})_i
	\mathop{\mathrm{d}\Omega} + 
	\int_{\Gamma} \pi_{N_h} \big( \widetilde{g_i} \big)
	\mathop{\mathrm{d}\Gamma} &= 0
	\quad \forall i \in \{1 : n - 1\},
	\label{eq:transient_compat_conds_disc_conc} \\
	\int_{\Gamma} \pi_{N_h} \big(
	\widetilde{g_J} \big) \mathop{\mathrm{d}\Gamma} &= 0.
	\label{eq:transient_compat_conds_disc_pot}
\end{align}
\end{subequations}
If $g_J$ depends on $\Phi_Z$ by a Butler--Volmer-type BC
\cref{eq:bv_linear}, then \cref{eq:transient_compat_conds_disc_pot}
constrains additive shifts in $\Phi_{Z,h}$.
Otherwise, assuming $g_J$ does not depend on $p$, $\bm{x}_{\nu}$ or $\Phi_Z$
and none of the $g_I$ depend on $\Phi_Z$,
the (consequently undetermined) additive constant in 
$\Phi_{Z,h}$ can be fixed by an extra constraint
$\int_{\Omega} \Phi_{Z,h} \mathop{\mathrm{d}\Omega} = 0$.
This leaves us to consider \cref{eq:transient_compat_conds_disc_conc}.

To study \cref{eq:transient_compat_conds_disc_conc} we assume a confined flow,
in the sense that the total normal fluxes
over $\Gamma$ are zero. In other words, we assume that the $g_i$ satisfy
\begin{equation} \label{eq:confined_assumption}
	\int_{\Gamma} \pi_{N_h} \big( \widetilde{g_i} \big)
	\mathop{\mathrm{d}\Gamma} = 0
	\quad \forall i \in \{1 : n - 1\},
\end{equation}
so that the compatibility conditions in 
\cref{eq:transient_compat_conds_disc_conc} become
\begin{equation} \label{eq:moles_conserved}
	\frac{\mathrm{d}}{\mathrm{d} t}
	\int_{\Omega} (\widetilde{\bm{c}_{\nu,h}})_i
	\mathop{\mathrm{d}\Omega} = 0
	\quad \forall i \in \{1 : n - 1\},
\end{equation}
which physically states that the total number of moles of all components is 
conserved.
Note that \cref{eq:confined_assumption} is satisfied whenever 
$g_i = \alpha_i g_J / F$ for some constants $\{\alpha_i\}_{i=1}^{n-1}$
(e.g.~in symmetric cells or lithium-ion battery cells),
since $g_J$ satisfies \cref{eq:transient_compat_conds_disc_pot}.
We now investigate how many of the $n-1$ conditions in \cref{eq:moles_conserved}
are actually independent.
Note that a general volumetric EOS relates
the total concentration $c_T$ to the partial molar volumes
$\bm{V}_{\nu}^{\top}$ and mole fractions $\bm{x}_{\nu}$ by
\cite{guggenheim1967thermodynamics}
\begin{equation} \label{eq:general_eos}
	c_T^{-1} = 
	\bm{V}_{\nu}^{\top} \bm{x}_{\nu}.
\end{equation}
Hence, multiplying \cref{eq:general_eos} by $c_T$, we see that the entries of 
$\widetilde{\bm{c}_{\nu,h}}$ are related by
\begin{equation} \label{eq:conc_constraints}
	1 = \widetilde{\bm{V}_{\nu}^{\top}} \widetilde{\bm{c}_{\nu,h}}.
\end{equation}
In light of \cref{eq:conc_constraints}, it is clear that the total molar
conservation conditions \cref{eq:moles_conserved}
may not all be independent, or may even be overdetermined.
This leads us to consider three pertinent cases, depending on the functional
dependence of $\bm{V}_{\nu}$ on $p$ and $\bm{x}_{\nu}$.

\textit{Case (i): The partial molar volumes are constant.}
In liquid mixtures it is often reasonable to approximate the partial molar volumes 
$\bm{V}_{\nu}$ as being constant \cite{druet2021global}.
It is then clear from \cref{eq:conc_constraints} that only
$n-2$ of the conditions in \cref{eq:moles_conserved} are independent.
If any $n-2$ of these conditions hold, then the final one will hold as well.
Thus, in this case, we must impose two additional constraints
(recall that the other constraint comes from taking $q_h = 1$ in
\cref{eq:disc_transient_ensosm_mavg}).
A reasonable choice is to enforce 
\cref{eq:mf_constraint} and a constraint that fixes the pressure such as
$\int_{\Omega} (p_h - \bar{p}) \mathop{\mathrm{d}\Omega}=0$.
If the partial molar volumes are constant then
$\bm{X}_{\nu}^0, \rho$ and $\bm{\psi}_Z$ are independent of $p$.
If $\eta, \zeta$ and $\bm{M}_Z$ are also modelled as being independent of $p$,
and $p$ does not appear in the BCs,
then the user-chosen value of $\bar{p}$ will not affect
the physical content of the model.

\textit{Case (ii): The partial molar volumes depend on
$\bm{x}_{\nu}$ but not $p$.}
If the partial molar volumes are non-constant 
(which happens only if they depend on $p$ or $\bm{x}_{\nu}$),
then the conditions in \cref{eq:moles_conserved} are generally independent.
One must then consider whether they can all be satisfied simultaneously.
If the partial molar volumes $\bm{V}_{\nu}$ depend on $\bm{x}_{\nu}$ only,
then the set of admissible concentrations $\bm{c}_{\nu}$ is parametrized by
$\bm{x}_{\nu}$.
Since only $n-2$ entries of $\bm{x}_{\nu}$ are independent, it does not seem
reasonable to expect that all $n-1$ conditions in 
\cref{eq:moles_conserved} can be satisfied simultaneously,
and we hypothesize that in this case the model problem \cref{eq:final_ensosm}
may not be well-posed.
We give an explicit example of how this case can produce an ill-posed model in 
\cref{sec:problematic_model}.

Challengingly, in real electrolytes the partial molar volumes can vary appreciably 
with composition but negligibly with pressure
\cite{newman2021electrochemical,wang2021potentiometric,wang2020shifting}.
To overcome the challenge of an ill-posed model one 
can explicitly include pressure dependence in $\bm{V}_{\nu}$.
However, if this dependence is too small, then to satisfy all $n-1$
conditions in \cref{eq:moles_conserved} the pressure may be forced to grow
unphysically large.
We hypothesize that this challenge could be overcome by allowing the volume
of $\Omega$ to vary in time
(i.e.~the total volume the electrolyte occupies may vary).
However, attempting this lies outside the scope of the present work,
since it would (very challengingly)
require both extending the model and numerically solving a moving boundary problem.
A simpler fix may be to use BCs that do not satisfy \cref{eq:confined_assumption}
and instead allow electrolyte to ``leak'' into or out of $\Omega$.

\textit{Case (iii): The partial molar volumes depend on
both $p$ and $\bm{x}_{\nu}$.}
As alluded to case (ii), if $\bm{V}_{\nu}$ depends on both
$p$ and $\bm{x}_{\nu}$, then the set of admissible concentrations 
$\bm{c}_{\nu}$ is parametrized by the $n-1$ independent quantities
$(p, \bm{x}_{\nu})$.
We therefore heuristically expect all $n-1$ conditions in
\cref{eq:moles_conserved} to be satisfiable.
In this case, only one additional constraint is needed
(from taking $q_h = 1$ in \cref{eq:disc_transient_ensosm_mavg}),
and we suggest employing \cref{eq:mf_constraint}.

Our analysis of cases (i)-(iii) only involves the total molar conservation
conditions \cref{eq:moles_conserved} and the volumetric EOS
\cref{eq:general_eos}.
In particular, our considerations are applicable to any
transient multicomponent fluid model with a general volumetric EOS.
We are not aware of these considerations being discussed elsewhere in the 
literature; particularly that of case (ii) where a seemingly physical choice
of EOS can lead to an ill-posed problem.
Although we are aware of numerical works
\cite{brunk2025structure,donev2015low}, which employ the EOS
\cref{eq:general_eos} with constant partial molar volumes,
they do not touch on the possibility of case (ii) and its associated 
challenges. Therefore, we believe that our considerations here may be of general 
interest to the multicomponent fluids community.

\subsection{Time-stepping methods} \label{sec:time_stepping}

The semi-discrete problem \cref{eq:disc_transient_ensosm},
together with strongly enforced BCs \cref{eq:disc_steady_bcs}
and appropriate integral constraints of \cref{sec:transient_constraints},
amounts to a nonlinear system of differential-algebraic equations (DAEs).
We choose to solve these DAEs using implicit Runge-Kutta methods, 
because of their desirable stability properties and ability to provide high-order 
temporal accuracy \cite{wanner1996solving}.

\section{Numerical examples} \label{sec:numerical}

Our subsequent numerical examples are implemented using Firedrake 
\cite{FiredrakeUserManual}.
We use Irksome \cite{farrell2021irksome,kirby2024extending} for time-stepping and
ngsPETSc \cite{ngspetsc,schoberl1997netgen} for mesh generation.
We solve the nonlinear systems with Newton's method \cite{deuflhard2011newton}
and the sparse direct solver MUMPS \cite{amestoy2001fully}.
Code can be found at
{\color{blue}\url{https://bitbucket.org/abaierr/multicomponent_electrolyte_code}}.

\subsection{Hull cell electroplating} \label{sec:hull_cell}

To demonstrate our methods in a physically realizable setting, we simulate
transient electroplating of a non-ideal binary electrolyte in a Hull cell
geometry in two and three spatial dimensions.
The two-dimensional domain $\Omega^{\textrm{2D}}$ is a right
trapezoid with vertices $(0, 0)$, $(0, 5)$, $(5, 5)$ and $(10, 0)$.
Here, a unit of length corresponds to a physical length of
$1 \si{\milli\metre}$.
We partition the boundary as
$\partial \Omega^{\textrm{2D}} = \Gamma_{p}^{\textrm{2D}} \cup 
\Gamma_{n}^{\textrm{2D}} \cup \Gamma_{w}^{\textrm{2D}}$
where $\Gamma_{p}^{\textrm{2D}}$ is the line segment between
$(0, 0)$ and $(0, 5)$ and denotes the positive electrode,
$\Gamma_{n}^{\textrm{2D}}$ is the line segment between
$(5, 5)$ and $(10, 0)$ and denotes the negative electrode,
and $\Gamma_{w}^{\textrm{2D}} = \partial \Omega^{\textrm{2D}} \setminus
(\Gamma_{p}^{\textrm{2D}} \cup \Gamma_{n}^{\textrm{2D}})$
are insulating walls.
The three-dimensional domain is the extrusion of the two-dimensional domain
by 5 units in the $z$-direction, i.e.~%
$\Omega^{\textrm{3D}} = \Omega^{\textrm{2D}} \times (0, 5)$
with $\Gamma_{p}^{\textrm{3D}} = \Gamma_{p}^{\textrm{2D}} \times (0, 5)$,
$\Gamma_{n}^{\textrm{3D}} = \Gamma_{n}^{\textrm{2D}} \times (0, 5)$
and $\Gamma_{w}^{\textrm{3D}} = \partial \Omega^{\textrm{3D}} \setminus
(\Gamma_{p}^{\textrm{3D}} \cup \Gamma_{n}^{\textrm{3D}})$.
In what follows we may omit superscripts
$\cdot^{\textrm{2D}}$ and $\cdot^{\textrm{3D}}$

The electrolyte is composed of ethyl-methyl-carbonate (\ce{EMC}) solvent and
lithium hexafluorophosphate (\ce{LiPF6}) salt, a mixture used in lithium-ion batteries.
In the notation of \cref{sec:cts}, there are $n=3$ species
(\ce{EMC}, \ce{Li+} and \ce{PF6-}) with molar masses
$\bm{m} = [104.105, 6.935, 144.97]^{\top} \si{\g\per\mole}$ and
equivalent charges
$\bm{z}=[0, 1, -1]^{\top}$.
We use the salt-charge transformation matrix \cref{eq:sc_trans_mat}
with $\bm{\nu}_1^{\top} = [1, 0, 0]$ and $\bm{\nu}_2^{\top} = [0, 1, 1]$.
This corresponds to neutralizing reactions
\ce{EMC <=> EMC} and \ce{Li+ + PF6- <=> LiPF6}.
Hence, under the salt-charge transformation, the mole fractions 
$(\bm{x}_{\nu})_i$, chemical potentials $(\bm{\mu}_{\nu})_i$, 
fluxes $(\bm{N}_{\nu})_i$ and so on,
represent those of \ce{EMC} for $i=1$ and \ce{LiPF6} for $i=2$.
We accordingly write $x_{\ce{EMC}} := (\bm{x}_{\nu})_1$,
$\mu_{\ce{EMC}} := (\bm{\mu}_{\nu})_1$,
$x_{\ce{LiPF6}} := (\bm{x}_{\nu})_2$,
$\mu_{\ce{LiPF6}} := (\bm{\mu}_{\nu})_2$
and so on for $(\bm{N}_{\nu})_i$.

The material properties are encoded in the functional dependence
of $\eta$, $\zeta$, $\rho$, $\bm{X}_{\nu}^0$ and $\bm{M}_Z$ on
$(T, p, \bm{x}_{\nu})$. We take an ambient temperature of
$T = 298.15\si{\kelvin}$. For the viscosities, we let
$\eta$ be a function of $\bm{x}_{\nu}$ as reported in
\cite[Table 2]{wang2020shifting}, and $\zeta = 10^{-6} \si{\pascal \second}$.
For $\rho$, $\bm{X}_{\nu}^0$ and $\bm{M}_Z$ we use the functional dependencies
from \cite[Table 1]{van2023structural}%
\footnote{The formula for $\kappa$ in \cite[Table 1]{van2023structural} has a typo
and in the notation of that paper should read
$\kappa = (48.93 y_e^{3/2} - 284.8 y_e^{5/2} + 817.8 y_e^4)^2$.
One can use the formulae from \cite{van2023structural} to construct
$\bm{M}_Z$, but note that the right side of
\cite[eq. 16.24]{van2023structural} is incorrect by a factor of $-1$.}.
This leads to $\rho$ depending on $\bm{x}_{\nu}$ only, and likewise for the
(non-constant) partial molar volumes, placing us in case (ii) of
\cref{sec:transient_constraints}.
The transport matrix $\bm{M}_Z$ also depends on $\bm{x}_{\nu}$ only, while
$\bm{X}_{\nu}^0$ depends on both $\bm{x}_{\nu}$ and $p$, with the (negligibly
small) dependence on $p$ arising due to the non-constant partial molar volumes.

For boundary conditions on $\Gamma := \partial \Omega$
we use \cref{eq:disc_steady_bcs_v} with
$g_{v_{\parallel}} = 0$.
For the normal fluxes
\cref{eq:disc_steady_bcs_N,eq:disc_steady_bcs_J},
we employ nonlinear Butler--Volmer BCs
\cite{dickinson2020butler,newman2021electrochemical}.
Note that $\mu_{\ce{LiPF6}}$ can be expressed as a known function of 
$(p, \bm{x}_{\nu})$ only, by analytically integrating $\bm{X}_{\nu}^0$ and 
$\bm{V}_{\nu}$ (cf.~\cref{eq:vt_ensosm_gradmu}).
The quantity $\Phi_Z - 0.5 \mu_{\ce{LiPF6}}/F$ then represents the potential of a
reference electrode that reacts with the electrolyte through
\ce{Li <=> Li+ + e-} \cite{van2023structural}.
We assume that the positive and negative electrodes
$\Gamma_{e} \textrm{ for } e \in \{p, n \}$ undergo the same reaction.
Letting $V_e$ denote the applied potential on electrode $e \in \{p, n \}$,
and $i_0(\bm{x}_{\nu})$ the exchange current density, we consider the BCs%
\begin{subequations}
\allowdisplaybreaks
\begin{alignat}{2}
	g_J &= - 2 i_0(\bm{x}_{\nu})
	\sinh \Bigg\{\frac{F 
	\big[ V_e - (\Phi_Z - 0.5 \mu_{\ce{LiPF6}} / F) \big]}{RT} 
	\Bigg\} \quad
	&&\textrm{ on } \Gamma_{e} \textrm{ for } e \in \{p, n \},
	\label{eq:bv_nonlinear_gJa} \\
	g_J &= 0 &&\textrm{ on } \Gamma_{w}, \label{eq:bv_nonlinear_gJb} \\
	g_2 &= g_J / (2 F) &&\textrm{ on } \Gamma. \label{eq:bv_nonlinear_g2}
\end{alignat}
\end{subequations}
Condition \cref{eq:bv_nonlinear_gJa} is a standard Butler--Volmer BC,
while \cref{eq:bv_nonlinear_gJb} expresses the insulating property of the walls.
Since $\bm{N} = \bm{Z}^{\top} \bm{N}_Z$, condition \cref{eq:bv_nonlinear_g2}
states that the normal flux of \ce{PF6-} is zero on $\Gamma$
(since the electrodes only react with \ce{Li+}).
We model $i_0$ by $i_0(\bm{x}_{\nu}) = i_0^{\ominus}
(x_{\ce{LiPF6}} / x_{\ce{LiPF6}}^{\ominus})^{1/2}$, to resemble commonly used
functional forms \cite{newman2021electrochemical}.
Moreover, we take $V_p = 0.1 \si{\volt}$, $V_n = 0 \si{\volt}$,
$i_0^{\ominus} = 10^4 \si{\ampere \meter^{-2}}$
and $x_{\ce{LiPF6}}^{\ominus} = 0.075$.

Since we are in case (ii) of \cref{sec:transient_constraints}
(non-constant partial molar volumes that do not depend on $p$),
we do not expect a well-posed problem if no-flux BCs are imposed on 
\ce{EMC} (i.e.~$g_1 = 0$ on $\Gamma$ in \cref{eq:disc_steady_bcs_N}).
Instead we impose BCs that allow some \ce{EMC} to ``leak'' from the
positive electrode. As a simple model for this, we assume a quadratic flux
profile of \ce{EMC} on the positive electrode,
with the magnitude of the flux an unknown to be solved for.
In particular, we take
\begin{subequations} \label{eq:emc_leak_bc}
\allowdisplaybreaks
\begin{alignat}{2}
	g_1 &= 0 &&\textrm{ on } \Gamma_{n} \cup \Gamma_{w}, \\
	g_1 &= \lambda_{\textrm{leak}} \cdot q_p \quad &&\textrm{ on } \Gamma_{p},
\end{alignat}
\end{subequations}
where $q_p : \Gamma_{p} \to \R$ is the unique quadratic function on 
$\Gamma_{p}$ that vanishes at its endpoints and is one at its midpoint, and
$\lambda_{\textrm{leak}} \in \R$ is a (time-dependent)
Lagrange multiplier that determines the amount of leakage.
The extra degree of freedom $\lambda_{\textrm{leak}}$ allows us to enforce
\cref{eq:mf_constraint} while also fixing the pressure mean, for which we take
$\int_{\Omega} p_h \mathop{\mathrm{d}\Omega}=0$.

As initial conditions we take a spatially uniform composition
$x_{\ce{LiPF6}}|_{t=0} = x_{\ce{LiPF6}}^{\ominus}$ and
$x_{\ce{EMC}}|_{t=0} = 1 - 2 x_{\ce{LiPF6}}^{\ominus}$
(cf.~\cref{eq:en_mfs_normalization}).
For compositions in this regime the Reynolds and Péclet numbers for the problem
are roughly $\textrm{Re} = 3 \times 10^{-5}$ and
$\textrm{Pe} = 9 \times 10^{-2}$.
We take all other unknowns $v$, $p$, $\bm{N}_{\nu}$, $J$, $\Phi_Z$
to be zero at $t=0$.
We numerically integrate \cref{eq:disc_transient_ensosm} up to a final time of
$172800$ seconds (two days) using the RadauIIA implicit Runge--Kutta 
method \cite{wanner1996solving} with two stages in the two-dimensional case
and one stage in the three-dimensional case, and with 200 timesteps.

We spatially discretize \cref{eq:disc_transient_ensosm} in
a non-dimensionalized form with augmentation parameter $\gamma = 10^{-2}$.
We employ a two-dimensional mesh of $5.7 \times 10^3$ triangles
with maximum cell diameter $h=0.125$, and a three-dimensional mesh of
$6.3 \times 10^3$ tetrahedra with maximum cell diameter $h=0.5$.
We expect singularities in the solution at corners of $\Omega$;
in the two-dimensional case we use a finer local cell diameter of 
$h_{\textrm{c}}=0.0125$ at the four vertices of $\Omega^{\textrm{2D}}$.
We employ the degree $k$ Taylor--Hood pair 
\cite{boffi1997three,taylor1973numerical} for $(V_h, P_h)$
and the $\mathbb{RT}_k$--$\mathbb{DG}_{k-1}$
\cite{raviart1977mixed} pair for $(N_h, X_h)$
with $k=4$ in two dimensions and $k=2$ in three dimensions.
The nonlinear systems at each timestep are solved using Newton's method
with an absolute tolerance on the residual of $10^{-10}$ in the Euclidean norm.
These systems consist of $1.3 \times 10^6$ unknowns in two dimensions
and $3.9 \times 10^5$ unknowns in three dimensions.
The greatest number of Newton iterations was taken at the first timestep
(10 iterations in two and three dimensions).
After the 2nd timestep, all Newton solves required at most 3 iterations.

\begin{figure}[h] \label{fig:hull_cell_2d}
	\caption{Streamlines of the \ce{EMC} flux $N_{\ce{EMC}}$
	(in black) and current density $J$ (in white) from the two-dimensional
	simulation of \cref{sec:hull_cell}. The domain is colored by
	the magnitude of $J$.
	}
	\centering
	\vspace{0.25em}
	\includegraphics[width=0.8\textwidth]{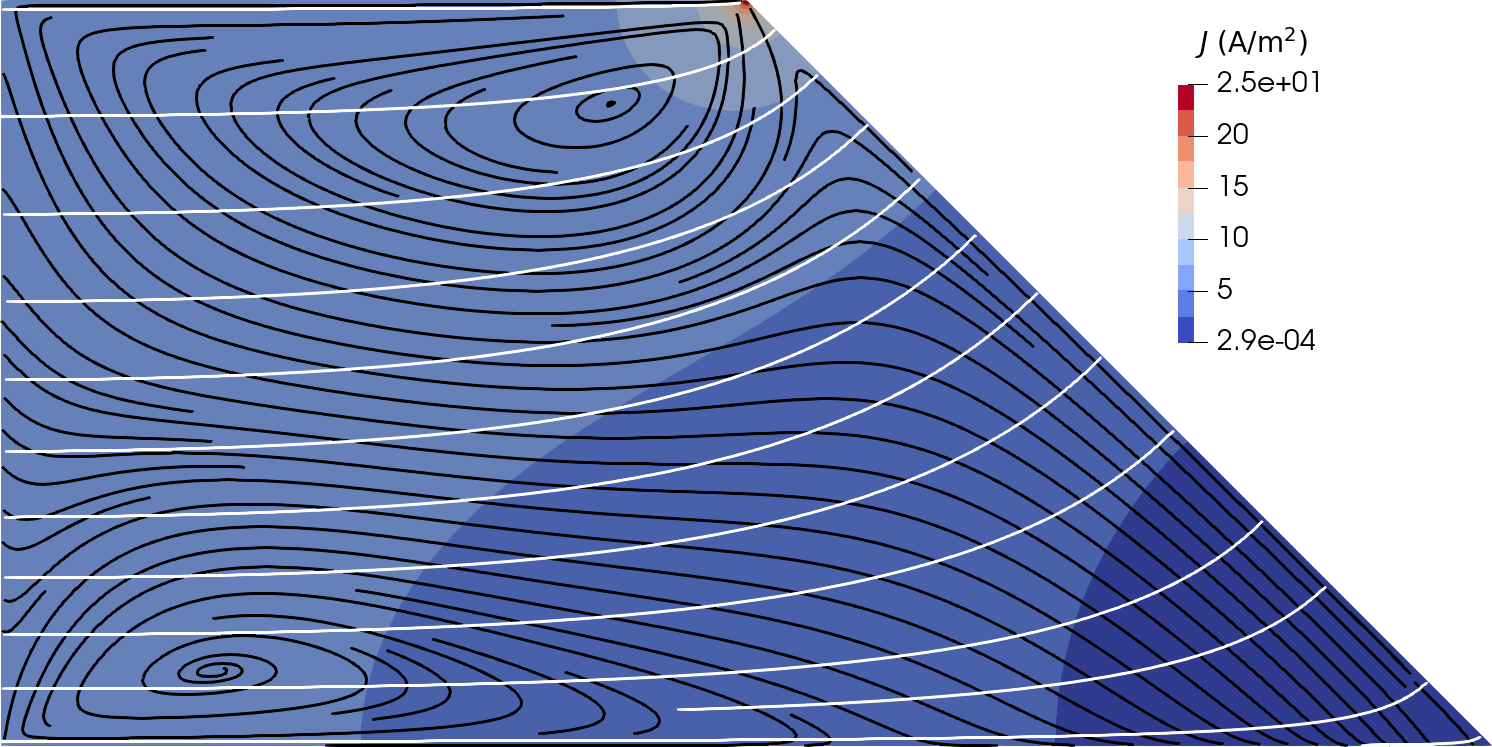}
\end{figure}

\begin{figure}[h] \label{fig:hull_cell_3d}
	\caption{\ce{EMC} flux $N_{\ce{EMC}}$ streamlines
		(colored by its magnitude)
		from the three-dimensional simulation of \cref{sec:hull_cell}.
		The domain walls are colored by the salt-charge potential $\Phi_Z$.
	}
	\centering
	\vspace{0.25em}
	\includegraphics[width=1.0\textwidth]{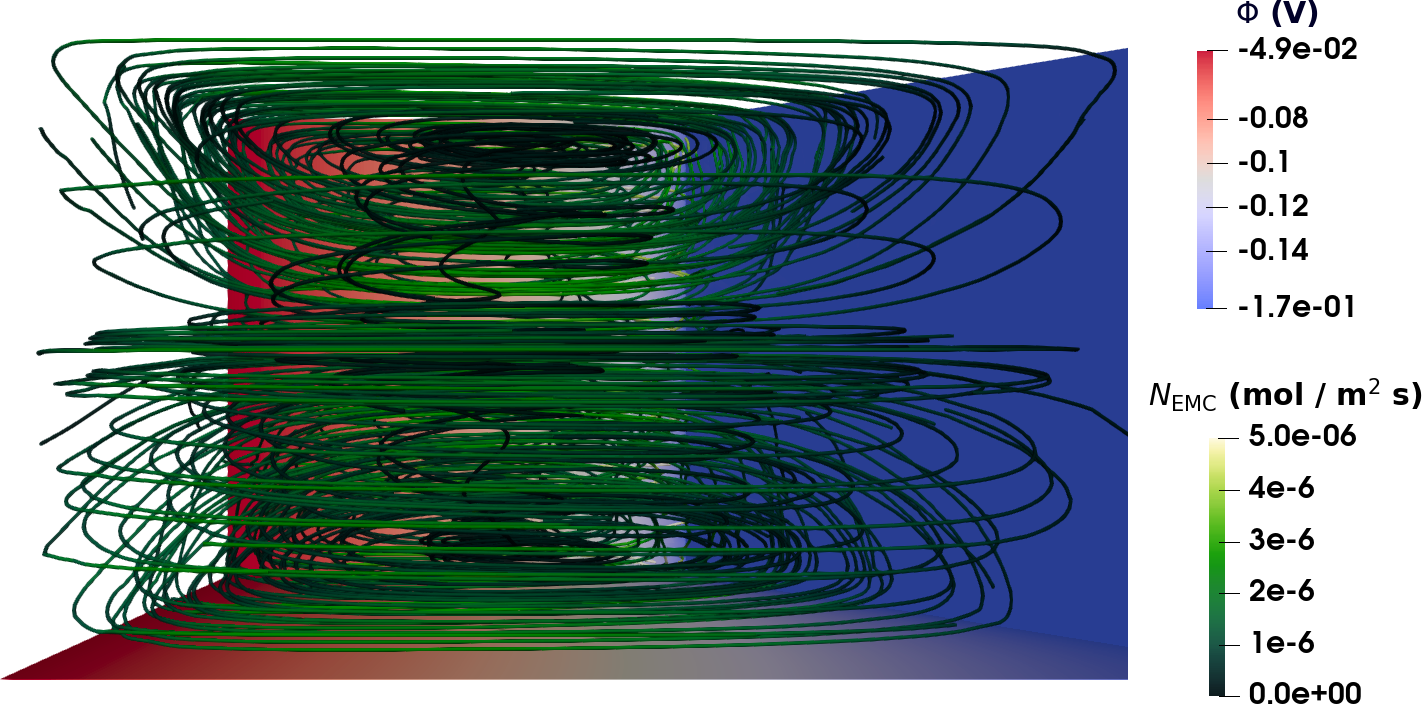}
\end{figure}

In \cref{fig:hull_cell_2d}, we plot streamlines of the \ce{EMC} flux 
$N_{\ce{EMC}}$ and current density
$J$ at time $t=64800 \si{\second}$ in two dimensions.
The current density expectedly flows from the positive to negative electrode
and appears to become singular at the top-right cell corner.
The two-dimensional plots also reveal two convective rolls forming in the
\ce{EMC} flux profile.
These convective rolls can also be seen in \cref{fig:hull_cell_3d}, where we plot 
streamlines of the \ce{EMC} flux at the final time $t=172800 \unit{\second}$ in 
three dimensions.

Owing to the leakage BCs \cref{eq:emc_leak_bc}, the
total number of moles of $\ce{EMC}$ (i.e.~the value of
$\int_{\Omega} c_1 \mathop{\mathrm{d}\Omega}$) varied by about
$0.08 \si{\percent}$ over the simulation time, in both two and three dimensions.
We also report maximum $L^2$-errors in the (nondimensionalized) mass-average 
constraint and mole fraction constraint:
\begin{align*}
	\mathcal{E}_1 &:= \max_{0 \leq t \leq T_{\textrm{final}}} \norm{v_h - 
	\widetilde{\bm{\psi}_Z^{\top}} \bm{N}_{Z,h}}_{L^2(\Omega)^d}
	\quad \textrm{ and } \quad
	\mathcal{E}_2 &:= \max_{0 \leq t \leq T_{\textrm{final}}}
	\norm{1 - \bm{\nu}_Z^{\top} \bm{x}_{\nu,h}}_{L^2(\Omega)}.
\end{align*}
We obtained values of
$(\mathcal{E}_1, \mathcal{E}_2) = (2.7 \times 10^{-4}, 5.6 \times 10^{-9})$
in two dimensions and
$(\mathcal{E}_1, \mathcal{E}_2) = (2.4 \times 10^{-2}, 2.5 \times 10^{-5})$
in three dimensions.

\subsection{Microfluidic rotating disk electrode} \label{sec:rotating_electrode}

We again consider the binary \ce{EMC}:\ce{LiPF6} electrolyte from 
\cref{sec:hull_cell}, but in a setting that showcases the flexibility of our 
method in applying different BCs.
We employ the same salt-charge transformation and material parameters
$\eta$, $\zeta$, $\rho$, $\bm{X}_{\nu}^0$ and $\bm{M}_Z$ at ambient temperature
$T = 298.15\si{\kelvin}$ as in \cref{sec:hull_cell}.
We consider a three-dimensional domain $\Omega$ representing a microfluidic
box containing a rotating disk electrode. Specifically, we take
$\Omega = \Omega_{\textrm{box}} \setminus \Omega_{\textrm{disk}}$ where
$\Omega_{\textrm{box}} = (-5, 5)^3$ and
$\Omega_{\textrm{disk}} = 
\{ (x, y, z) \in \R^3 : x^2 + y^2 \leq 1, -0.05 \leq z \leq 0.05 \}$.
A unit of length corresponds to a physical length of
$12.5 \si{\micro\metre}$.
We decompose $\Gamma := \partial \Omega$ as
$\Gamma = \Gamma_{p} \cup \Gamma_{n} \cup \Gamma_{w}$ where
$\Gamma_{p} = \{ (x, y, z) \in \partial \Omega_{\textrm{box}} : 
z = \pm 5 \}$ denotes the top and bottom walls of the box,
$\Gamma_{w} = \partial \Omega_{\textrm{box}} \setminus \Gamma_{p}$
the four side walls of the box, and
$\Gamma_{n} = \partial \Omega_{\textrm{disk}}$ the disk boundary.
The surfaces $\Gamma_{p}$, $\Gamma_{n}$, $\Gamma_{w}$ represent, respectively,
the positive electrode, negative electrode, and insulating walls.

We solve the steady problem \cref{eq:disc_steady_ensosm} with the
following BCs. We first impose that the disk rotates with fixed angular frequency
$\dot{\theta} \in \R$; this manifests in our model through the boundary
condition \cref{eq:disc_steady_bcs_v} on $v$ with $g_{v_{\parallel}}$ given by
\begin{subequations}
\allowdisplaybreaks
\begin{alignat}{2}
	g_{v_{\parallel}} &= \dot{\theta} (-y, x, 0) \quad && \textrm{ for }
	(x, y, z) \in \partial \Omega_{\textrm{disk}} = \Gamma_{n}, \\
	g_{v_{\parallel}} &= 0 &&\textrm{ on } \Gamma \setminus \partial 
	\Omega_{\textrm{disk}}.
\end{alignat}
\end{subequations}
We choose $\dot{\theta} = 28.79 \si{\second^{-1}}$, which leads to
Reynolds and Péclet numbers of roughly $\textrm{Re} = 3 \times 10^{-3}$ and
$\textrm{Pe} = 1 \times 10^{1}$.
For the $\ce{EMC}$ flux $N_{h,\ce{EMC}} := (\bm{N}_{\nu,h})_1$ we impose
a zero normal flux condition \cref{eq:disc_steady_bcs_N} with
$g_1 = 0$ on $\Gamma$.
For the current density flux, instead of enforcing \cref{eq:disc_steady_bcs_J}
for some function $g_J$, we enforce 
$J_h \cdot n_{\Gamma} = 2 F (N_{h,\ce{LiPF6}} \cdot n_{\Gamma})$
on $\Gamma$,
where $N_{h,\ce{LiPF6}} := (\bm{N}_{\nu,h})_2$
(note that we implement this as a strongly enforced BC on $J_h \cdot n_{\Gamma}$
and not $N_{h,\ce{LiPF6}} \cdot n_{\Gamma}$).
As in \cref{sec:hull_cell}, this BC enforces that the normal flux of 
\ce{PF6-} is zero on $\Gamma$.
For the $\ce{LiPF6}$ flux
we enforce a zero normal flux condition on the insulating walls
$N_{h,\ce{LiPF6}} \cdot n_{\Gamma} = 0$ on $\Gamma_{w}$.
However, instead of prescribing the value of
$N_{\ce{LiPF6}} \cdot n_{\Gamma}$ on the electrodes $\Gamma_p \cup \Gamma_n$,
we weakly enforce
\begin{equation} \label{eq:rotating_disc_bcs}
	x_{h,\ce{LiPF6}} = x_{\ce{LiPF6}}^{\ominus,e}
	\textrm{ on } \Gamma_e \textrm{ for } e \in \{p, n\},
\end{equation}
where $x_{\ce{LiPF6}}^{\ominus,p} = 0.082$ and
$x_{\ce{LiPF6}}^{\ominus,n} = 0.068$.
This is done by modifying 
\cref{eq:disc_steady_ensosm_osm} through
the addition of boundary terms, and by enlarging the space of
test functions $(\bm{W}_h)_2$ to those with zero normal trace on $\Gamma_{w}$ only.
Specifically, instead of \cref{eq:disc_steady_ensosm_osm}, we consider:
\begin{align} \label{eq:disc_steady_ensosm_osm_modified}
\allowdisplaybreaks
\begin{split}
	&\Bigg(
	\begin{bmatrix}
		\bm{x}_{\nu,h} \\
		\Phi_{Z,h}
	\end{bmatrix},
	\nabla \cdot
	\begin{bmatrix}
		RT \widetilde{\bm{X}_{\nu}^0} \bm{W}_h \\
		F \norm{\bm{z}} K_h^{\top}
	\end{bmatrix}
	\Bigg)_{\Omega}
	-
	\Bigg(p_h,
	\nabla \cdot
	\Bigg\{
	\widetilde{\bm{\psi}_{Z}^{\top}}
	\begin{bmatrix}
		\bm{W}_h \\
		K_h^{\top}
	\end{bmatrix}
	- \widetilde{\bm{V}_{\nu}^{\top}} \bm{W}_h
	\Bigg\}
	\Bigg)_{\Omega} \\
	&+ \ \Bigg(
	\gamma  \widetilde{\bm{\psi}_{Z}} v_h,
	\begin{bmatrix}
		\bm{W}_h \\
		K_h^{\top}
	\end{bmatrix}
	\Bigg)_{\Omega}
	=
	\Bigg(
	\widetilde{\bm{M}_Z^{\gamma}} \bm{N}_{Z,h},
	\begin{bmatrix}
		\bm{W}_h \\
		K_h^{\top}
	\end{bmatrix} \Bigg)_{\Omega} \\
	&+
	RT \cdot x_{\ce{LiPF6}}^{\ominus,p} 
	\Big\langle \widetilde{(\bm{X}_{\nu}^0)_{22}}
	, (\bm{W}_h)_2 \cdot n_{\Gamma} \Big\rangle_{\Gamma_p}
	+ RT \cdot x_{\ce{LiPF6}}^{\ominus,n}
	\Big\langle \widetilde{(\bm{X}_{\nu}^0)_{22}}
	, (\bm{W}_h)_2 \cdot n_{\Gamma} \Big\rangle_{\Gamma_n}
	\\
	&\quad \forall
	\begin{bmatrix}
		\bm{W}_h \\
		K_h^{\top}
	\end{bmatrix} \in (N_{0h} \times N_h \times N_{0h})
	\quad \textrm{with }
	(\bm{W}_h)_2 \cdot n_{\Gamma} = 0 \textrm{ on } \Gamma_w.
\end{split}
\end{align}
With the present material data,
$\bm{X}_{\nu}^0$ is a diagonal matrix.
One thus verifies
by taking $(\bm{W}_h)_2$ to be non-zero in 
\cref{eq:disc_steady_ensosm_osm_modified},
that for $e \in \{p, n\}$ the following is weakly enforced:
\begin{equation} \label{eq:rotating_disc_bc_true}
	RT \cdot x_{h,\ce{LiPF6}} \cdot \widetilde{(\bm{X}_{\nu}^0)_{22}}
	- \underbrace{p_h \cdot \Big\{ (\widetilde{\bm{\psi}_{Z}^{\top})_2} - 
	\widetilde{(\bm{V}_{\nu}^{\top})_2} \Big\}}_{:=I_p}
	= RT \cdot x_{\ce{LiPF6}}^{\ominus,e} \cdot \widetilde{(\bm{X}_{\nu}^0)_{22}}
	\quad \textrm{on } \Gamma_e.
\end{equation}
Neglecting $I_p$%
\footnote{Dimensional analysis reveals that $I_p$ is smaller than the other terms
in \cref{eq:rotating_disc_bc_true} by a factor of $10^{-8}$.}
in \cref{eq:rotating_disc_bc_true} and dividing
by $RT \cdot \widetilde{(\bm{X}_{\nu}^0)_{22}}$ then leads to
\cref{eq:rotating_disc_bcs} as desired.

We discretize \cref{eq:disc_steady_ensosm}
with aforementioned BCs and modification \cref{eq:disc_steady_ensosm_osm_modified}
in a non-dimensionalized form with augmentation parameter $\gamma = 10^{-2}$.
As constraints we impose \cref{eq:mf_constraint} along with
$\int_{\Omega} p_h \mathop{\mathrm{d}\Omega}=0$.
We employ a curved mesh of degree $4$ with
$1.8 \times 10^4$ tetrahedra and maximum local cell diameter of $h=0.1$
on the disk boundary.
We employ the degree $k=4$ Taylor--Hood pair 
\cite{boffi1997three,taylor1973numerical} for $(V_h, P_h)$
and the $\mathbb{RT}_k$--$\mathbb{DG}_{k-1}$
\cite{raviart1977mixed} pair for $(N_h, X_h)$.
The nonlinear system consists of $1.1 \times 10^6$ unknowns and was solved using 
Newton's method with an absolute tolerance on the residual of $10^{-10}$ in the 
Euclidean norm.
We first applied Newton's method on a coarse discretization
(with a non-curved mesh and order $k=2$ spaces),
where as an initial guess we set
$x_{\ce{LiPF6}} = 0.075$
and $x_{\ce{EMC}} = 1 - 2 x_{\ce{LiPF6}}$
(cf.~\cref{eq:en_mfs_normalization}) and we set all other unknowns to be zero.
Convergence was reached in 6 iterations. We then used the coarse solution
as an initial guess for the fine discretization
(i.e.~with a curved mesh and degree $k=4$ spaces, as above),
for which Newton's method converged in 3 iterations.

\begin{figure}[h] \label{fig:rotating_electrode}
\caption{Streamlines of the \ce{LiPF6} flux $N_{\ce{LiPF6}}$
	(colored by its magnitude)
	and current density $J$ (colored in a transparent yellow)
	above $\Omega_{\textrm{disk}}$
	for the numerical experiment of \cref{sec:rotating_electrode}.
	The disk is also colored by the magnitude of $N_{\ce{LiPF6}}$.
}
\centering
\vspace{0.25em}
\includegraphics[width=1.0\textwidth]{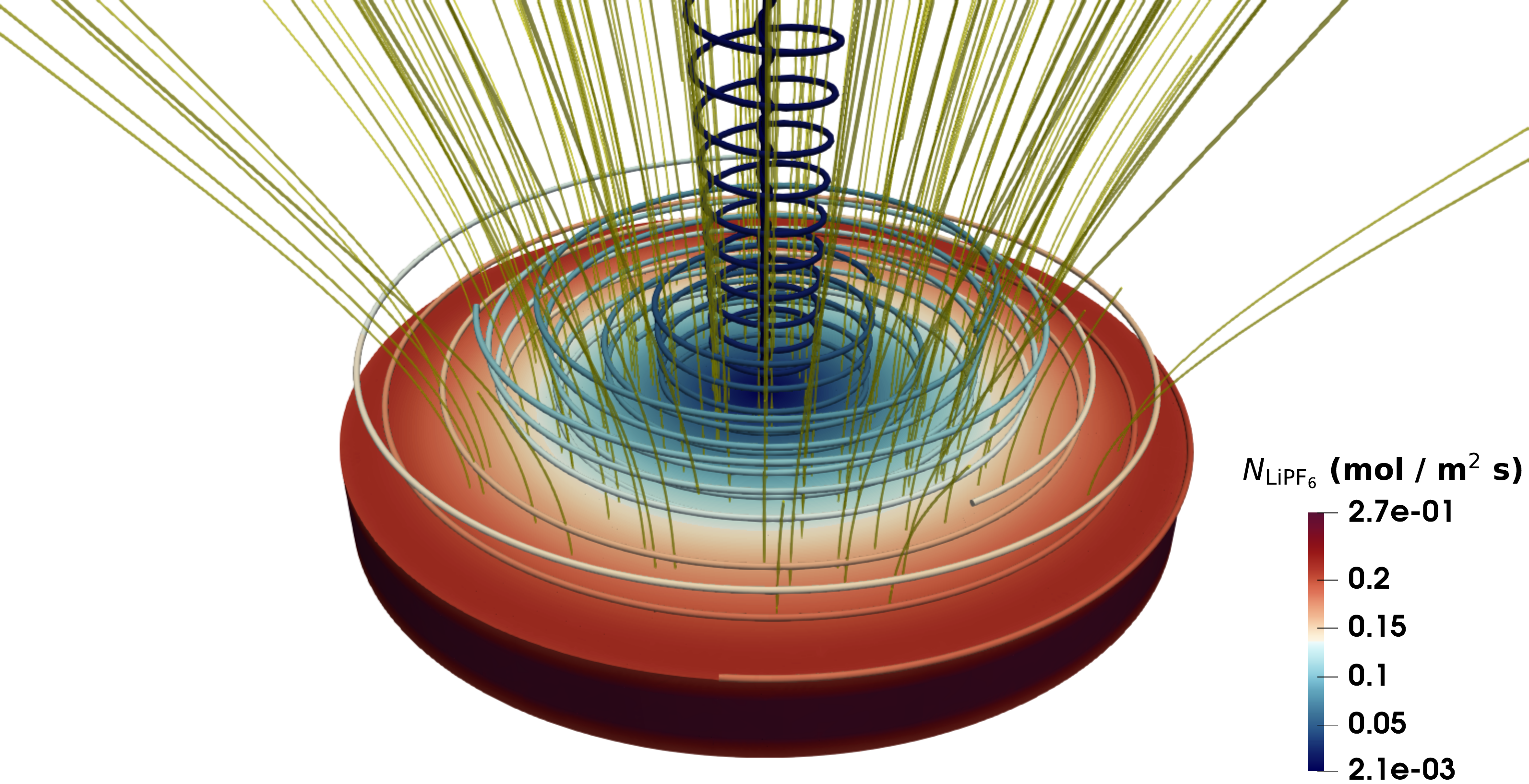}
\end{figure}

In \cref{fig:rotating_electrode} we plot streamlines of the $\ce{LiPF6}$ flux and 
current density $J$.
The rotation of the disk induces a swirling flow of \ce{LiPF6}, while the
current density flows smoothly from the positive electrode to the negative
disk electrode. These observations qualitatively reflect the expected
behaviour of the solution as induced by our choice of BCs.
Moreover, the $L^2$-errors in the (nondimensionalized) mass-average 
constraint, mole fraction constraint, and weakly enforced BCs
\cref{eq:rotating_disc_bcs} are
\begin{align*}
	\norm{v_h - \widetilde{\bm{\psi}_Z^{\top}} \bm{N}_{Z,h}}_{L^2(\Omega)^d}
	&= 1.3 \times 10^{-3}, \\
	\norm{1 - \bm{\nu}_Z^{\top} \bm{x}_{\nu,h}}_{L^2(\Omega)}
	&= 3.7 \times 10^{-7}, \\
	\norm{x_{h,\ce{LiPF6}} - x_{\ce{LiPF6}}^{\ominus,e}}
	_{L^2(\Gamma_p \cup \Gamma_n)}
	&=5.0 \times 10^{-5}.
\end{align*}
In particular, the small error in the weakly enforced BC
\cref{eq:rotating_disc_bcs} affirms the ability of our algorithm to 
handle a combination of flux, Dirichlet, and tangential velocity BCs.

\subsection{Cosolvent imbalances} \label{sec:cosolvent}

An appealing property of our numerical method is that it accounts
for the multicomponent nature of the electrolyte solvent
(i.e.~the neutral species in which the salts are dissolved),
the effects of which may be important for battery modelling 
\cite{jung2023overpotential}.
To demonstrate this, we 
consider an electrolyte comprised of \ce{LiPF6} salt dissolved in 
two solvents, namely ethyl-methyl-carbonate (\ce{EMC}) and
ethylene-carbonate (\ce{EC}).
The resulting mixture has $n=4$ species
(\ce{EMC}, \ce{EC}, \ce{Li+} and \ce{PF6-}) with molar masses
$\bm{m} = [104.105, 88.062, 6.935, 144.97]^{\top} \si{\g\per\mole}$ and
equivalent charges
$\bm{z}=[0, 0, 1, -1]^{\top}$.
We use the salt-charge transformation matrix \cref{eq:sc_trans_mat}
with $\bm{\nu}_1^{\top} = [1, 0, 0, 0]$,
$\bm{\nu}_2^{\top} = [0, 1, 0, 0]$
and $\bm{\nu}_3^{\top} = [0, 0, 1, 1]$.
This corresponds to neutralizing reactions
\ce{EMC <=> EMC}, \ce{EC <=> EC} and \ce{Li+ + PF6- <=> LiPF6}.
Hence, under the salt-charge transformation, the mole fractions 
$(\bm{x}_{\nu})_i$, chemical potentials $(\bm{\mu}_{\nu})_i$, 
fluxes $(\bm{N}_{\nu})_i$ and so on,
represent those of \ce{EMC} for $i=1$, \ce{EC} for $i=2$ and \ce{LiPF6} for $i=3$.
\begin{figure}[h] \label{fig:cosolvent_streamlines}
	\caption{Streamlines of the \ce{EMC} flux $N_{\ce{EMC}}$
		(in black) and \ce{EC} flux $N_{\ce{EC}}$ (in white) from the 
		simulation of \cref{sec:cosolvent}. The domain is colored by
		the shear viscosity $\eta$.
	}
	\centering
	\vspace{0.25em}
	\includegraphics[width=0.8\textwidth]{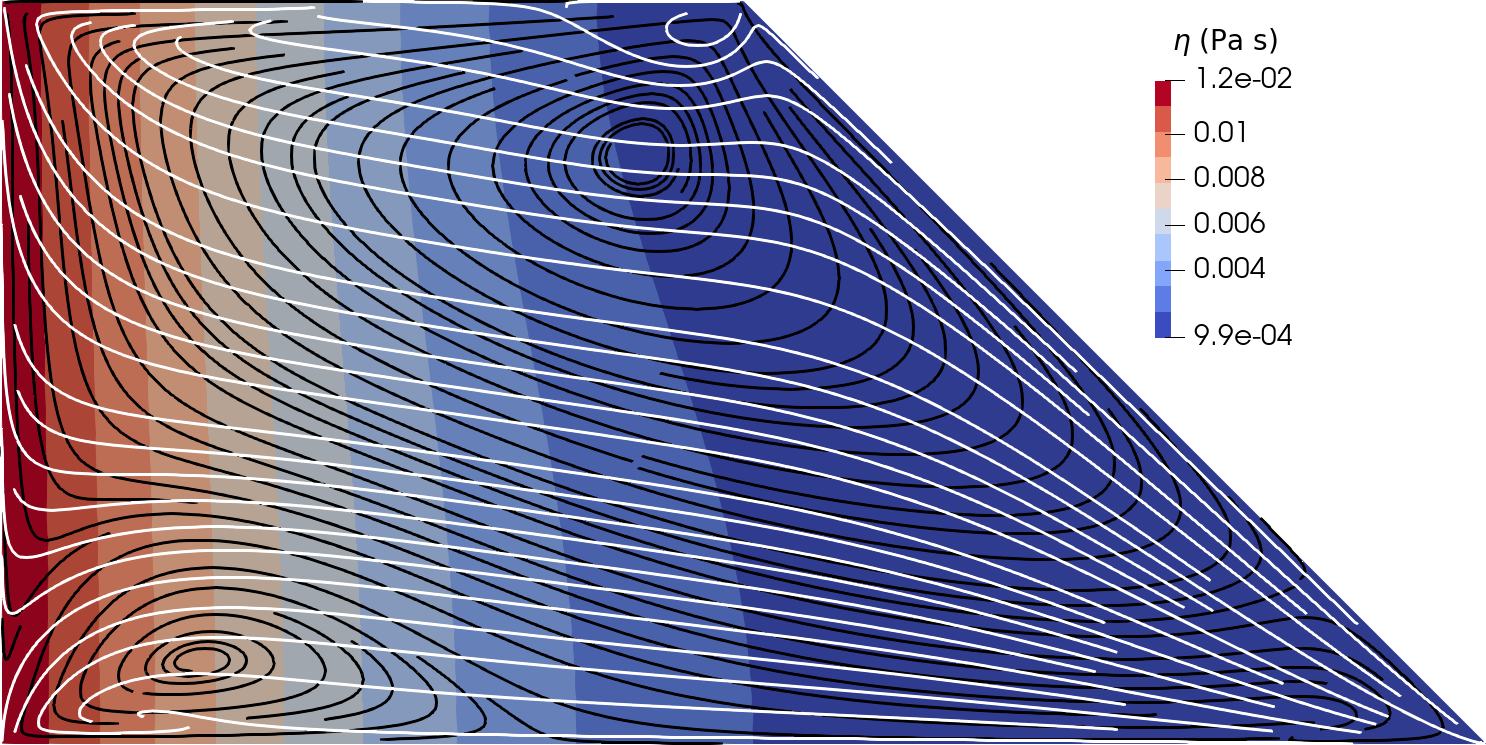}
\end{figure}
We take an ambient temperature of
$T = 298.15\si{\kelvin}$.
We let $\rho$ be a function of $\bm{x}_{\nu}$ as reported in
\cite{rungta2025quantifying}.
An expression for $\eta$ as a function of $\bm{x}_{\nu}$ was obtained by
fitting viscosity data reported in \cite{rungta2025quantifying} to a degree
7 bivariate polynomial in the variables
$\sqrt{x_{\ce{EMC}} / (x_{\ce{EMC}} + x_{\ce{EC}})}$,
$\sqrt{x_{\ce{EC}} / (x_{\ce{EMC}} + x_{\ce{EC}})}$.
We further take $\zeta = 10^{-6} \si{\pascal \second}$, and
$\bm{X}_{\nu}^0$ is obtained by assuming an ideal mixture.
We compute $\bm{M}_Z$ by assuming constant
Stefan--Maxwell diffusivities, the values of which are obtained
from the supplementary information of \cite{phelan2024applying} for
$1 \si{\mole \per \litre}$ of $\ce{LiPF6}$ in vol:vol 7:3 \ce{EMC}:\ce{EC} solvent.

We take $\Omega$ to be the same two-dimensional Hull cell domain
as in \cref{sec:hull_cell}.
For BCs we use \cref{eq:disc_steady_bcs_v} with
$g_{v_{\parallel}} = 0$.
For \cref{eq:disc_steady_bcs_N,eq:disc_steady_bcs_J}
we use linearized Butler--Volmer BCs and we neglect the
$\mu_{\ce{LiPF6}}$ term in the overpotential (cf.~\cref{eq:bv_linear}):
\begin{subequations}
\allowdisplaybreaks
	\begin{alignat}{2}
		g_J &= -i_0 F (V_{\textrm{e}} - \Phi_Z) / (RT)
		\quad &&\textrm{ on } \Gamma_{e} \textrm{ for } e \in \{p, n \},
		 \\
		g_J &= 0 &&\textrm{ on } \Gamma_{w}, \\
		g_3 &= g_J / (2 F) &&\textrm{ on } \Gamma,
	\end{alignat}
\end{subequations}
with $V_p = 0.1 \si{\volt}$, $V_n = 0 \si{\volt}$ and
$i_0 = 10^4 \si{\ampere \meter^{-2}}$.
In \cref{eq:disc_steady_bcs_N} we further take $g_2 = 0$ on $\Gamma$,
i.e.~no normal flux of \ce{EC}, while for \ce{EMC} we impose the
leak BCs \cref{eq:emc_leak_bc}.

We numerically solve the transient problem \cref{eq:disc_transient_ensosm}.
As initial conditions we take a spatially uniform composition
$x_{\ce{LiPF6}}|_{t=0} = 0.077$, 
$x_{\ce{EMC}}|_{t=0} = 0.509$ and
$x_{\ce{EC}}|_{t=0} = 1 - 0.509 - (2 \cdot 0.077) = 0.337$
(cf.~\cref{eq:en_mfs_normalization}).
In this regime the Reynolds and Péclet numbers for the problem
are roughly $\textrm{Re} = 3 \times 10^{-5}$ and
$\textrm{Pe} = 9 \times 10^{-1}$.
We take all other unknowns $v$, $p$, $\bm{N}_{\nu}$, $J$, $\Phi_Z$
to be zero at $t=0$.
As constraints we impose \cref{eq:mf_constraint} along with
$\int_{\Omega} p_h \mathop{\mathrm{d}\Omega}=0$.
We use the same mesh, finite 
element spaces $(V_h, P_h)$ and $(N_h, X_h)$, and time-stepping scheme as in the 
two-dimensional case of \cref{sec:hull_cell}. 
The nonlinear systems at each timestep consist of $1.7 \times 10^6$ unknowns and 
are solved using Newton's method with an absolute tolerance on the residual of 
$10^{-10}$ in the Euclidean norm.
The first two timesteps required at most 5 Newton iterations, and all subsequent 
timesteps required at most 3 iterations. Note that the initial Newton solves in 
\cref{sec:hull_cell} required more iterations due to the nonlinear BCs 
\cref{eq:bv_nonlinear_gJa}.

In \cref{fig:cosolvent_streamlines} 
we plot streamlines of the \ce{EMC} and \ce{EC} fluxes
$N_{\ce{EMC}}$, $N_{\ce{EC}}$ at time $t = 64800$.
The flux profiles are emphatically different, underscoring the multicomponent 
nature of the solvent.
Moreover, in \cref{fig:cosolvent_streamlines} we color the domain by
the shear viscosity $\eta = \eta(\bm{x}_{\nu})$, which changes by an order of
magnitude across the cell due to the spatially varying $\ce{EMC}$:$\ce{EC}$ ratio.
Specifically, at $t=64800$ the ratio
$x_{\ce{EMC}} / x_{\ce{EC}} \approx 1.4$ near the positive (left) electrode and
$x_{\ce{EMC}} / x_{\ce{EC}} \approx 1.6$ near the negative (right) electrode.

\section{Conclusions} \label{sec:conc}
We have presented a broad family of finite element algorithms for numerically
solving the electroneutral NSOSM equations.
To the best of our knowledge, this is the first paper in the finite element
literature on electroneutral NSOSM flow.
The flexibility of our algorithms in handling transient and steady flow under
different boundary conditions was substantiated in our numerical experiments.
Our numerical experiment involving \ce{EMC}-\ce{EC}-\ce{LiPF6} flow
also demonstrated the scientific potential of our algorithms in studying how the 
multicomponent nature of electrolytes may impact, for example, locally varying 
material properties of the mixture.

\appendix
\section{A problematic model} \label{sec:problematic_model}
Consider a binary mixture on a fixed domain $\Omega$
with molar concentrations $c_1, c_2$, total concentration 
$c_T = c_1 + c_2$ and molar fractions $x_i := c_i / c_T$ with $x_1 + x_2 = 1$.
Let $x := x_1$.
We assume a volumetric EOS
\begin{equation} \label{eq:illposed_eos}
	c_T = A + B x,
\end{equation}
with $A$ and $B$ non-zero constants satisfying $A > 0$ and $A + B > 0$ so that 
$c_T > 0$.
This seemingly benign EOS assumes the total concentration
is linear in composition; an experimentalist might reasonably make such an 
approximation by fitting experimental data.
One can verify that the (non-constant) partial molar volumes are
\begin{equation*}
	V_1 = \frac{A + B(2 x - 1)}{(A + Bx)^2}, \quad
	V_2 = \frac{A + 2 B x}{(A + Bx)^2} \quad
	\textrm{and} \quad
	\frac{1}{c_T} = x_1 V_1 + x_2 V_2.
\end{equation*}
Assume that the total number of moles of each species is conserved, so that
\begin{equation} \label{eq:illposed_c}
	\frac{\mathrm{d}}{\mathrm{d} t}
	\int_{\Omega} c_i
	\mathop{\mathrm{d}\Omega} = 0
	\quad \forall i \in \{1, 2\}.
\end{equation}
Since $c_T = c_1 + c_2$, we can use \cref{eq:illposed_c}
and use \cref{eq:illposed_eos} to deduce that
\begin{equation} \label{eq:illposed_x}
	\frac{\mathrm{d}}{\mathrm{d} t}
	\int_{\Omega} x
	\mathop{\mathrm{d}\Omega} = 0.
\end{equation}
Moreover, since $c_1 = x c_T$ we can use
\cref{eq:illposed_c,eq:illposed_eos,eq:illposed_x}
to deduce that
\begin{equation} \label{eq:illposed_x2}
	\frac{\mathrm{d}}{\mathrm{d} t}
	\int_{\Omega} x^2
	\mathop{\mathrm{d}\Omega} = 0.
\end{equation}
Let $\bar{x} := \int_{\Omega} x\mathop{\mathrm{d}\Omega} / 
\int_{\Omega} 1\mathop{\mathrm{d}\Omega}$ denote the spatial mean of $x$.
Since $\Omega$ is assumed not to depend on $t$,
\cref{eq:illposed_x} implies that $\bar{x}$ is a constant.
Note that \cref{eq:illposed_x,eq:illposed_x2} also imply
\begin{equation*}
\frac{\mathrm{d}}{\mathrm{d} t}
\int_{\Omega} (x - \bar{x})^2 \mathop{\mathrm{d}\Omega}
= \frac{\mathrm{d}}{\mathrm{d} t}\Big(
\int_{\Omega} x^2 \mathop{\mathrm{d}\Omega}
- \int_{\Omega} \bar{x}^2 \mathop{\mathrm{d}\Omega} \Big)
= 0.
\end{equation*}
Thus $\int_{\Omega} (x - \bar{x})^2 \mathop{\mathrm{d}\Omega} = C$ for 
some constant $C$.
Now, if $x$ is initially spatially uniform, then $C=0$. This then implies that
$x = \bar{x}$ a.e.~in $\Omega$ for all time, and $\bar{x}$ is a constant.

To conclude, any model which assumes EOS \cref{eq:illposed_eos} and
conservation properties \cref{eq:illposed_c},
with $\Omega$ independent of $t$, is seemingly either (i) 
ill-posed or (ii) well-posed, but 
incapable of making physically meaningful predictions, since a solution that 
initially has a spatially uniform composition retains that
uniform composition over all time.

\bibliographystyle{siamplain}

\begin{thebibliography}{10}

\bibitem{amestoy2001fully}
{\sc P.~R. Amestoy, I.~S. Duff, J.-Y. L'Excellent, and J.~Koster}, {\em A fully
  asynchronous multifrontal solver using distributed dynamic scheduling}, SIAM
  J. Matrix Anal. Appl., 23 (2001), pp.~15--41.

\bibitem{aznaran2024finite}
{\sc F.~R. Aznaran, P.~E. Farrell, C.~W. Monroe, and A.~J. Van-Brunt}, {\em
  Finite element methods for multicomponent convection-diffusion}, IMA J.
  Numer. Anal.,  (2024), p.~drae001.

\bibitem{baier2024high}
{\sc A.~Baier-Reinio and P.~E. Farrell}, {\em High-order finite element methods
  for three-dimensional multicomponent convection-diffusion}, arXiv preprint
  arXiv:2408.17390,  (2024).

\bibitem{balakrishnan2014fluctuating}
{\sc K.~Balakrishnan, A.~L. Garcia, A.~Donev, and J.~B. Bell}, {\em Fluctuating
  hydrodynamics of multispecies nonreactive mixtures}, Phys. Rev. E, 89 (2014),
  p.~013017.

\bibitem{bard2022electrochemical}
{\sc A.~J. Bard and L.~R. Faulkner}, {\em Electrochemical {M}ethods:
  {F}undamentals and {A}pplications}, John Wiley \& Sons, New York, 2nd~ed.,
  2001.

\bibitem{bartlett2008bioelectrochemistry}
{\sc P.~N. Bartlett}, {\em Bioelectrochemistry: {F}undamentals, {E}xperimental
  {T}echniques and {A}pplications}, John Wiley \& Sons, Chichester, 2008.

\bibitem{ngspetsc}
{\sc J.~Betteridge, P.~E. Farrell, M.~Hochsteger, C.~Lackner, J.~Schöberl,
  S.~Zampini, and U.~Zerbinati}, {\em {ngsPETSc}: A coupling between
  {NETGEN}/{NGSolve} and {PETSc}}, J. Open Source Softw., 9 (2024), p.~7359.

\bibitem{bhattacharjee2015fluctuating}
{\sc A.~K. Bhattacharjee, K.~Balakrishnan, A.~L. Garcia, J.~B. Bell, and
  A.~Donev}, {\em Fluctuating hydrodynamics of multi-species reactive
  mixtures}, J. Chem. Phys., 142 (2015), p.~224107.

\bibitem{bird2002transport}
{\sc R.~B. Bird, W.~E. Stewart, and E.~N. Lightfoot}, {\em Transport
  {P}henomena}, John Wiley \& Sons, 2nd~ed., 2002.

\bibitem{bizeray2016resolving}
{\sc A.~M. Bizeray, D.~A. Howey, and C.~W. Monroe}, {\em Resolving a
  discrepancy in diffusion potentials, with a case study for {L}i-ion
  batteries}, J. Electrochem. Soc., 163 (2016), p.~E223.

\bibitem{boettcher2020potentially}
{\sc S.~W. Boettcher, S.~Z. Oener, M.~C. Lonergan, Y.~Surendranath, S.~Ardo,
  C.~Brozek, and P.~A. Kempler}, {\em Potentially confusing: potentials in
  electrochemistry}, ACS Energy Lett., 6 (2020), pp.~261--266.

\bibitem{boffi1997three}
{\sc D.~Boffi}, {\em Three-dimensional finite element methods for the {S}tokes
  problem}, SIAM J. Numer. Anal., 34 (1997), pp.~664--670.

\bibitem{bortels1996multi}
{\sc L.~Bortels, J.~Deconinck, and B.~Van Den~Bossche}, {\em The
  multi-dimensional upwinding method as a new simulation tool for the analysis
  of multi-ion electrolytes controlled by diffusion, convection and migration.
  {P}art 1. {S}teady state analysis of a parallel plane flow channel}, J.
  Electroanal. Chem., 404 (1996), pp.~15--26.

\bibitem{braukhoff2022entropy}
{\sc M.~Braukhoff, I.~Perugia, and P.~Stocker}, {\em An entropy structure
  preserving space-time formulation for cross-diffusion systems: analysis and
  {G}alerkin discretization}, SIAM J. Numer. Anal., 60 (2022), pp.~364--395.

\bibitem{brezzi1985two}
{\sc F.~Brezzi, J.~Douglas, and L.~D. Marini}, {\em Two families of mixed
  finite elements for second order elliptic problems}, Numer. Math., 47 (1985),
  pp.~217--235.

\bibitem{brosa2022continuum}
{\sc F.~Brosa~Planella, W.~Ai, A.~M. Boyce, A.~Ghosh, I.~Korotkin, S.~Sahu,
  V.~Sulzer, R.~Timms, T.~G. Tranter, M.~Zyskin, et~al.}, {\em A continuum of
  physics-based lithium-ion battery models reviewed}, Prog. Energy, 4 (2022),
  p.~042003.

\bibitem{brunk2025structure}
{\sc A.~Brunk, A.~J{\"u}ngel, and M.~Luk{\'a}{\v{c}}ov{\'a}-Medvid'ov{\'a}},
  {\em A structure-preserving numerical method for quasi-incompressible
  {Navier--Stokes--Maxwell--Stefan} systems}, arXiv preprint arXiv:2504.11892,
  (2025).

\bibitem{burman2003bunsen}
{\sc E.~Burman, A.~Ern, and V.~Giovangigli}, {\em Bunsen flame simulation by
  finite elements on adaptively refined, unstructured triangulations}, Combust.
  Theory Model., 8 (2003), p.~65.

\bibitem{butler1924studies}
{\sc J.~A.~V. Butler}, {\em Studies in heterogeneous equilibria. {P}art
  {II}.—{T}he kinetic interpretation of the {N}ernst theory of electromotive
  force}, Trans. Faraday Soc., 19 (1924), pp.~729--733.

\bibitem{carnes2008local}
{\sc B.~Carnes and G.~F. Carey}, {\em Local boundary value problems for the
  error in {FE} approximation of non-linear diffusion systems}, Internat. J.
  Numer. Methods Engrg., 73 (2008), pp.~665--684.

\bibitem{correa2023new}
{\sc C.~I. Correa, G.~N. Gatica, and R.~Ruiz-Baier}, {\em New mixed finite
  element methods for the coupled {S}tokes and {Poisson--Nernst--Planck}
  equations in {B}anach spaces}, ESAIM Math. Model. Numer. Anal., 57 (2023),
  pp.~1511--1551.

\bibitem{de2013non}
{\sc S.~R. De~Groot and P.~Mazur}, {\em Non-{E}quilibrium {T}hermodynamics},
  Dover Publications, Inc., New York, 1984.

\bibitem{deuflhard2011newton}
{\sc P.~Deuflhard}, {\em Newton {M}ethods for {N}onlinear {P}roblems}, Springer
  Science \& Business Media, Berlin, Heidelberg, 2011.

\bibitem{dickinson2020butler}
{\sc E.~J. Dickinson and A.~J. Wain}, {\em The {Butler-Volmer} equation in
  electrochemical theory: {O}rigins, value, and practical application}, J.
  Electroanal. Chem., 872 (2020), p.~114145.

\bibitem{donev2015low}
{\sc A.~Donev, A.~Nonaka, A.~K. Bhattacharjee, A.~L. Garcia, and J.~B. Bell},
  {\em Low {M}ach number fluctuating hydrodynamics of multispecies liquid
  mixtures}, Phys. Fluids, 27 (2015), p.~037103.

\bibitem{donev2014low}
{\sc A.~Donev, A.~Nonaka, Y.~Sun, T.~Fai, A.~Garcia, and J.~Bell}, {\em Low
  {M}ach number fluctuating hydrodynamics of diffusively mixing fluids},
  Commun. Appl. Math. Comput. Sci., 9 (2014), pp.~47--105.

\bibitem{donev2019fluctuating}
{\sc A.~Donev, A.~J. Nonaka, C.~Kim, A.~L. Garcia, and J.~B. Bell}, {\em
  Fluctuating hydrodynamics of electrolytes at electroneutral scales}, Phys.
  Rev. Fluids, 4 (2019), p.~043701.

\bibitem{dreyer2013overcoming}
{\sc W.~Dreyer, C.~Guhlke, and R.~M{\"u}ller}, {\em Overcoming the shortcomings
  of the {N}ernst--{P}lanck model}, Phys. Chem. Chem. Phys., 15 (2013),
  pp.~7075--7086.

\bibitem{druet2021global}
{\sc P.-{\'E}. Druet}, {\em Global--in--time existence for liquid mixtures
  subject to a generalised incompressibility constraint}, J. Math. Anal. Appl.,
  499 (2021), p.~125059.

\bibitem{ellingsrud2025splitting}
{\sc A.~J. Ellingsrud, P.~Benedusi, and M.~Kuchta}, {\em A splitting,
  discontinuous {G}alerkin solver for the cell-by-cell electroneutral
  {N}ernst--{P}lanck framework}, SIAM J. Sci. Comput., 47 (2025),
  pp.~B477--B504.

\bibitem{erdey1930theorie}
{\sc T.~Erdey-Gr{\'u}z and M.~Volmer}, {\em Zur {T}heorie der {W}asserstoff
  {{\"U}}berspannung}, Z. Phys. Chem., 150A (1930), pp.~203--213.

\bibitem{ern1994multicomponent}
{\sc A.~Ern and V.~Giovangigli}, {\em Multicomponent {T}ransport {A}lgorithms},
  vol.~24, Springer Berlin, Heidelberg, 1994.

\bibitem{ern1998thermal}
{\sc A.~Ern and V.~Giovangigli}, {\em Thermal diffusion effects in hydrogen-air
  and methane-air flames}, Combust. Theory Model., 2 (1998), p.~349.

\bibitem{ern2017finite}
{\sc A.~Ern and J.-L. Guermond}, {\em Finite element quasi-interpolation and
  best approximation}, ESAIM Math. Model. Numer. Anal., 51 (2017),
  pp.~1367--1385.

\bibitem{ern2021finiteI}
{\sc A.~Ern and J.-L. Guermond}, {\em Finite {E}lements {I}: {A}pproximation
  and {I}nterpolation}, Springer, Cham, Switzerland, 2021.

\bibitem{ern2021finiteII}
{\sc A.~Ern and J.-L. Guermond}, {\em Finite {E}lements {II}: {G}alerkin
  {A}pproximation, {E}lliptic and {M}ixed {PDE}s}, Springer, Cham, Switzerland,
  2021.

\bibitem{farrell2021irksome}
{\sc P.~E. Farrell, R.~C. Kirby, and J.~Marchena-Menendez}, {\em {I}rksome:
  {A}utomating {Runge--Kutta} time-stepping for finite element methods}, ACM
  Trans. Math. Softw., 47 (2021), pp.~1--26.

\bibitem{feireisl2016mathematical}
{\sc E.~Feireisl, D.~Hilhorst, H.~Petzeltov{\'a}, and P.~Tak{\'a}{\v{c}}}, {\em
  Mathematical analysis of variable density flows in porous media}, J. Evol.
  Equ., 16 (2016), pp.~1--19.

\bibitem{fick1855ueber}
{\sc A.~Fick}, {\em {\"U}ber {D}iffusion}, Annalen der Physik, 170 (1855),
  pp.~59--86.

\bibitem{giovangigli2012multicomponent}
{\sc V.~Giovangigli}, {\em Multicomponent {F}low {M}odeling}, Birkhäuser,
  Boston, 1999.

\bibitem{doble2007perry}
{\sc D.~W. Green and R.~H. Perry}, {\em Perry’s {C}hemical {E}ngineers’
  {H}andbook}, McGraw Hill Professional, 8th~ed., 2007.

\bibitem{guggenheim1967thermodynamics}
{\sc E.~A. Guggenheim}, {\em Thermodynamics: {A}n {A}dvanced {T}reatment for
  {C}hemists and {P}hysicists}, North-Holland Books, Amsterdam, 5th~ed., 1967.

\bibitem{FiredrakeUserManual}
{\sc D.~A. Ham and {26 Others}}, {\em Firedrake {U}ser {M}anual}, 1st~ed.,
  2023, \url{https://doi.org/10.25561/104839}.

\bibitem{helfand1960inversion}
{\sc E.~Helfand}, {\em On inversion of the linear laws of irreversible
  thermodynamics}, J. Chem. Phys., 33 (1960), pp.~319--322.

\bibitem{jung2023overpotential}
{\sc T.~Jung, A.~A. Wang, and C.~W. Monroe}, {\em Overpotential from cosolvent
  imbalance in battery electrolytes: {LiPF6 in EMC: EC}}, ACS omega, 8 (2023),
  pp.~21133--21144.

\bibitem{jungel2015boundedness}
{\sc A.~J{\"u}ngel}, {\em The boundedness-by-entropy method for cross-diffusion
  systems}, Nonlinearity, 28 (2015), p.~1963.

\bibitem{jungel2019convergence}
{\sc A.~J{\"u}ngel and O.~Leingang}, {\em Convergence of an implicit {E}uler
  {G}alerkin scheme for {P}oisson--{M}axwell--{S}tefan systems}, Adv. Comput.
  Math., 45 (2019), pp.~1469--1498.

\bibitem{kirby2024extending}
{\sc R.~C. Kirby and S.~P. MacLachlan}, {\em Extending {I}rksome: improvements
  in automated {Runge--Kutta} time stepping for finite element methods}, ACM
  Trans. Math. Softw.,  (2024).

\bibitem{kraaijeveld1995modelling}
{\sc G.~Kraaijeveld, V.~Sumberova, S.~Kuindersma, and H.~Wesselingh}, {\em
  Modelling electrodialysis using the {Maxwell--Stefan} description}, Chem.
  Eng. J., 57 (1995), pp.~163--176.

\bibitem{kraaijeveld1993negative}
{\sc G.~Kraaijeveld and J.~A. Wesselingh}, {\em Negative {M}axwell--{S}tefan
  diffusion coefficients}, Ind. Eng. Chem. Res., 32 (1993), pp.~738--742.

\bibitem{krishna2015uphill}
{\sc R.~Krishna}, {\em Uphill diffusion in multicomponent mixtures}, Chem. Soc.
  Rev., 44 (2015), pp.~2812--2836.

\bibitem{krishna2019diffusing}
{\sc R.~Krishna}, {\em Diffusing uphill with {J}ames {C}lerk {M}axwell and
  {J}osef {S}tefan}, Chem. Eng. Sci., 195 (2019), pp.~851--880.

\bibitem{krishna1997maxwell}
{\sc R.~Krishna and J.~A. Wesselingh}, {\em The {M}axwell--{S}tefan approach to
  mass transfer}, Chem. Eng. Sci., 52 (1997), pp.~861--911.

\bibitem{longo2012finite}
{\sc A.~Longo, M.~Barsanti, A.~Cassioli, and P.~Papale}, {\em A finite element
  {G}alerkin/least-squares method for computation of multicomponent
  compressible--incompressible flows}, Comput. \& Fluids, 67 (2012),
  pp.~57--71.

\bibitem{maxwell1866}
{\sc J.~C. Maxwell}, {\em On the dynamical theory of gases}, Phil. Trans. R.
  Soc.,  (1866), pp.~49--88.

\bibitem{mcleod2014mixed}
{\sc M.~McLeod and Y.~Bourgault}, {\em Mixed finite element methods for
  addressing multi-species diffusion using the {M}axwell--{S}tefan equations},
  Comput. Meth. Appl. Mech. Eng., 279 (2014), pp.~515--535.

\bibitem{nedelec1986new}
{\sc J.-C. N{\'e}d{\'e}lec}, {\em A new family of mixed finite elements in
  {$\mathbb{R}^3$}}, Numer. Math., 50 (1986), pp.~57--81.

\bibitem{newman2021electrochemical}
{\sc J.~Newman and N.~P. Balsara}, {\em Electrochemical {S}ystems}, John Wiley
  \& Sons, Hoboken, NJ, 4th~ed., 2021.

\bibitem{newman1965mass}
{\sc J.~Newman, D.~Bennion, and C.~W. Tobias}, {\em Mass transfer in
  concentrated binary electrolytes}, Ber. Bunsenges. Phys. Chem., 69 (1965),
  pp.~608--612.

\bibitem{onsager1931reciprocal}
{\sc L.~Onsager}, {\em Reciprocal relations in irreversible processes. {I}.},
  Phys. Rev., 37 (1931), pp.~405--426.

\bibitem{onsager1931reciprocal2}
{\sc L.~Onsager}, {\em Reciprocal relations in irreversible processes. {II}.},
  Phys. Rev., 38 (1931), pp.~2265--2279.

\bibitem{peraud2016low}
{\sc J.-P. P{\'e}raud, A.~Nonaka, A.~Chaudhri, J.~B. Bell, A.~Donev, and A.~L.
  Garcia}, {\em Low {M}ach number fluctuating hydrodynamics for electrolytes},
  Phys. Rev. Fluids, 1 (2016), p.~074103.

\bibitem{pethica2007electrostatic}
{\sc B.~A. Pethica}, {\em Are electrostatic potentials between regions of
  different chemical composition measurable? the {Gibbs--Guggenheim} principle
  reconsidered, extended and its consequences revisited}, Phys. Chem. Chem.
  Phys., 9 (2007), pp.~6253--6262.

\bibitem{phelan2024applying}
{\sc C.~Phelan, J.~Swallow, and R.~Weatherup}, {\em Applying the
  {M}axwell-{S}tefan diffusion framework to multicomponent battery
  electrolytes}, ChemRxiv preprint,  (2024).

\bibitem{prohl2010convergent}
{\sc A.~Prohl and M.~Schmuck}, {\em Convergent finite element discretizationsof
  the {Navier--Stokes--Nernst--Planck--Poisson} system}, ESAIM Math. Model.
  Numer. Anal., 44 (2010), pp.~531--571.

\bibitem{raviart1977mixed}
{\sc P.-A. Raviart and J.-M. Thomas}, {\em A mixed finite element method for
  2-nd order elliptic problems}, in {M}athematical {A}spects of {F}inite
  {E}lement {M}ethods, vol.~606 of {L}ecture {N}otes in {M}ath., Springer,
  Berlin, 1977, pp.~{292--315}.

\bibitem{richardson2022charge}
{\sc G.~W. Richardson, J.~M. Foster, R.~Ranom, C.~P. Please, and A.~M. Ramos},
  {\em Charge transport modelling of {L}ithium-ion batteries}, European J.
  Appl. Math., 33 (2022), pp.~983--1031.

\bibitem{roy2023scalable}
{\sc T.~Roy, J.~Andrej, and V.~A. Beck}, {\em A scalable {DG} solver for the
  electroneutral {N}ernst--{P}lanck equations}, J. Comput. Phys., 475 (2023),
  p.~111859.

\bibitem{rungta2025quantifying}
{\sc R.~Rungta, P.~Slowikowski, A.~Gardner, D.~Persa, and C.~W. Monroe}, {\em
  Quantifying volumetric expansion and bulk moduli of carbonate cosolvents with
  lithium salts}, in preparation.

\bibitem{schoberl1997netgen}
{\sc J.~Sch\"oberl}, {\em {NETGEN}: An advancing front {2D}/{3D}-mesh generator
  based on abstract rules}, Computing and Visualization in Science, 1 (1997),
  pp.~41--52.

\bibitem{scott1985norm}
{\sc L.~R. Scott and M.~Vogelius}, {\em Norm estimates for a maximal right
  inverse of the divergence operator in spaces of piecewise polynomials}, ESAIM
  Math. Model. Numer. Anal., 19 (1985), pp.~111--143.

\bibitem{sijabat2019maxwell}
{\sc R.~Sijabat, M.~De~Groot, S.~Moshtarikhah, and J.~Van Der~Schaaf}, {\em
  {Maxwell--Stefan} model of multicomponent ion transport inside a monolayer
  {Nafion} membrane for intensified chlor-alkali electrolysis}, J. Appl.
  Electrochem., 49 (2019), pp.~353--368.

\bibitem{srivastava2023staggered}
{\sc I.~Srivastava, D.~R. Ladiges, A.~J. Nonaka, A.~L. Garcia, and J.~B. Bell},
  {\em Staggered scheme for the compressible fluctuating hydrodynamics of
  multispecies fluid mixtures}, Phys. Rev. E, 107 (2023), p.~015305.

\bibitem{stefan1871gleichgewicht}
{\sc J.~Stefan}, {\em {\"U}ber das {G}leichgewicht und die {B}ewegung,
  insbesondere die {D}iffusion von {G}asgemengen}, Sitzber. Akad. Wiss. Wien.,
  63 (1871), pp.~63--124.

\bibitem{sun2019entropy}
{\sc Z.~Sun, J.~A. Carrillo, and C.-W. Shu}, {\em An entropy stable high-order
  discontinuous {G}alerkin method for cross-diffusion gradient flow systems},
  Kinet. Relat. Models, 12 (2019), pp.~885--908.

\bibitem{taylor1973numerical}
{\sc C.~Taylor and P.~Hood}, {\em A numerical solution of the
  {N}avier--{S}tokes equations using the finite element technique}, Comput. \&
  Fluids, 1 (1973), pp.~73--100.

\bibitem{van2022augmented}
{\sc A.~Van-Brunt, P.~E. Farrell, and C.~W. Monroe}, {\em Augmented
  saddle-point formulation of the steady-state {S}tefan--{M}axwell diffusion
  problem}, IMA J. Numer. Anal., 42 (2022), pp.~3272--3305.

\bibitem{van2022consolidated}
{\sc A.~Van-Brunt, P.~E. Farrell, and C.~W. Monroe}, {\em Consolidated theory
  of fluid thermodiffusion}, AlChE J., 68 (2022), p.~e17599.

\bibitem{van2023structural}
{\sc A.~Van-Brunt, P.~E. Farrell, and C.~W. Monroe}, {\em Structural
  electroneutrality in {Onsager--Stefan--Maxwell} transport with charged
  species}, Electrochim. Acta, 441 (2023), p.~141769.

\bibitem{vanag2009cross}
{\sc V.~K. Vanag and I.~R. Epstein}, {\em Cross-diffusion and pattern formation
  in reaction--diffusion systems}, Phys. Chem. Chem. Phys., 11 (2009),
  pp.~897--912.

\bibitem{wang2021potentiometric}
{\sc A.~A. Wang, A.~B. Gunnarsd{\'o}ttir, J.~Fawdon, M.~Pasta, C.~P. Grey, and
  C.~W. Monroe}, {\em Potentiometric {MRI} of a superconcentrated lithium
  electrolyte: testing the irreversible thermodynamics approach}, ACS Energy
  Lett., 6 (2021), pp.~3086--3095.

\bibitem{wang2020shifting}
{\sc A.~A. Wang, T.~Hou, M.~Karanjavala, and C.~W. Monroe}, {\em
  Shifting-reference concentration cells to refine composition-dependent
  transport characterization of binary lithium-ion electrolytes}, Electrochim.
  Acta, 358 (2020), p.~136688.

\bibitem{wanner1996solving}
{\sc G.~Wanner and E.~Hairer}, {\em Solving {O}rdinary {D}ifferential
  {E}quations {II}}, vol.~375, Springer Berlin Heidelberg, 1996.

\bibitem{weber2014critical}
{\sc A.~Z. Weber, R.~L. Borup, R.~M. Darling, P.~K. Das, T.~J. Dursch, W.~Gu,
  D.~Harvey, A.~Kusoglu, S.~Litster, M.~M. Mench, et~al.}, {\em A critical
  review of modeling transport phenomena in polymer-electrolyte fuel cells}, J.
  Electrochem. Soc., 161 (2014), p.~F1254.

\bibitem{weber2008mathematical}
{\sc A.~Z. Weber and C.~Delacourt}, {\em Mathematical modelling of cation
  contamination in a proton-exchange membrane}, Fuel Cells, 8 (2008),
  pp.~459--465.

\bibitem{wesselingh2000mass}
{\sc J.~Wesselingh and R.~Krishna}, {\em Mass {T}ransfer in {M}ulticomponent
  {M}ixtures}, Delft University Press, Delft, Netherlands, 2000.

\bibitem{xie2020effective}
{\sc D.~Xie and B.~Lu}, {\em An effective finite element iterative solver for a
  {Poisson--Nernst--Planck} ion channel model with periodic boundary
  conditions}, SIAM J. Sci. Comput., 42 (2020), pp.~B1490--B1516.

\end{thebibliography}

\end{document}